\documentclass{article}
\usepackage{amsmath}
\usepackage{amssymb}
\usepackage{mathtools}
\usepackage{verbatim}
\usepackage{enumerate}
\usepackage{graphicx}
\usepackage{float}
\usepackage{overpic}
\usepackage{subfigure}
\usepackage{caption}
\usepackage{relsize}
\usepackage{geometry}
\usepackage{multirow}
\usepackage{tikz}
\usepackage{epstopdf}

\newcommand\ii{\,\mathrm{i}}
\newcommand\bR{\mathbb{R}}

\newcommand\bZ{\mathbb{Z}}

\newcommand{\tdiag}{\text{diag}}

\begin{document}
\title{Localized Sparsifying Preconditioner for Periodic Indefinite Systems}
\author{Fei Liu$^\sharp$ and Lexing Ying$^{\dagger\sharp}$\\
  $\dagger$ Department of Mathematics, Stanford University\\
  $\sharp$ Institute for Computational and Mathematical Engineering, Stanford University
}
\date{}
\maketitle

\begin{abstract}
This paper introduces the localized sparsifying preconditioner for the
pseudospectral approximations of indefinite systems on periodic
structures. The work is built on top of the recently proposed
sparsifying preconditioner with two major modifications. First, the local
potential information is utilized to improve the accuracy of the
preconditioner. Second, an FFT based method to compute the local
stencil is proposed to reduce the setup time of the
algorithm. Numerical results show that the iteration number of this
improved method grows only mildly as the problem size grows, which
implies that solving pseudospectral approximation systems is
computationally as efficient as solving sparse systems, up to a mildly
growing factor.
\end{abstract}



\section{Introduction}
\label{sec:introuction}
This paper is concerned with the numerical solution of highly
indefinite systems on periodic structures with periodic boundary
condition
\begin{equation}
\label{eqn:model}
(-\Delta + v(x))u(x) = f(x),\quad x\in [0,1)^d,
\end{equation}
where $d$ is the dimension, $f(x)$ is the right-hand side, and $u(x)$
is the unknown. $v(x)$ is the potential that can take negative
values. In the case of the periodic Helmholtz equation, $v(x)$ is
$-(\omega/c(x))^2$ where $\omega$ is the angular frequency and $c(x)$
is the velocity field, while for the periodic Schr\"odinger equation,
$v(x)$ is a rescaling of $v_{\text{ext}}(x) - E$ where
$v_{\text{ext}}(x)$ is the external potential field and $E$ is the
energy shift.

Solving \eqref{eqn:model} numerically is a challenging task since the
system can be highly indefinite, which makes most of the classic
iterative solvers no longer effective. Moreover, the solution
typically has a highly oscillatory pattern and it requires large
number of unknowns for accurate numerical approximations due to the
Nyquist theorem.

The simplest way to solve \eqref{eqn:model} numerically is to adopt
the standard second order central difference scheme, which results a
sparse system, then the sparse direct methods, such as the nested
dissection method
\cite{george1973nested,duff1983multifrontal,liu1992multifrontal}, can
be applied directly. However, the dispersion relation given by the
standard second order central difference scheme is not accurate
enough, which leads to a poor approximation of the solution. One way
to fix this is to use higher order difference schemes. The problem is
that, higher order schemes require larger stencil supports, therefore
the effectiveness of the sparse direct solvers cannot be leveraged.

A more natural way to discretize \eqref{eqn:model} is to use the
pseudospectral method \cite{gottliebspectral,trefethenspectral} with
Fourier basis. The pseudospectral differentiation scheme requires only
a few points per oscillation of the solution to give an accurate
dispersion relation. However, the stencil induced by the scheme is
not local, thus the direct sparse solvers cannot be applied directly.

Recently in \cite{spspd,spspc}, the sparsifying preconditioners are
proposed to address the issue of balancing the accuracy and the
sparsity pattern of the discretized systems. The main idea is to
numerically convert the dense linear system derived from some accurate
scheme into a sparse system and use the inverse of the sparse system
as a preconditioner for the dense one. The numerical results in
\cite{spspc} show a satisfying iteration number for solving the
Lippmann-Schwinger equation. However in \cite{spspd}, the iteration
number needed to solve the indefinite system \eqref{eqn:model} is not
as small, because the periodic boundary condition implies a higher
requirement for the accuracy of the dispersion relation.

This paper is a follow-up work of \cite{spspd}. We propose the
localized sparsifying preconditioner which takes the local information
of the potential $v(x)$ into consideration in order to give a more
accurate sparse approximation of the non-sparse pseudospectral
system. In addition, an FFT based method for computing the local
stencils is proposed to accelerate the setup process of the
preconditioner.

The rest of the paper is organized as follows. We first formulate the
pseudospectral discrete system to be solved in Section
\ref{sec:formulation}. Section \ref{sec:review} briefly reviews the
previous work in \cite{spspd}. In Section \ref{sec:localized} we
present the modifications in this new work. Numerical results are
given in Section \ref{sec:numerical}.  Section \ref{sec:conclusion}
concludes with some extra discussions.

\section{Formulation of the Pseudospectral System}
\label{sec:formulation}
This section introduces the pseudospectral discretization for solving
\eqref{eqn:model}. Discretizing \eqref{eqn:model} with $n$ points
along each dimension results in a uniform Cartesian grid on $[0,1)^d$,
  which can be indexed by the set
\begin{align*}
J: = \{(j_1,\dots,j_d):0\le j_1,\dots,j_d < n\}.
\end{align*}
The corresponding Cartesian grid is denoted by $hJ$ where $h=1/n$ is
the step size. $h$ is chosen such that we have at least four points
per oscillation. We also introduce a grid in the Fourier domain
\begin{align*}
K := \{(k_1,\dots,k_d):-n/2\le k_1,\dots, k_d < n/2\}.
\end{align*}
Given array $z=\{z_j: j\in J\}$ defined on grid $J$ and $\hat{z} = \{\hat{z}_k: k \in K \}$ defined on grid $K$, the Fourier and inverse Fourier transforms $F$ and $F^{-1}$ are defined as
\begin{gather*}
\hat{z}_k = (Fz)_k = \dfrac{1}{n^{d/2}}\sum_{j\in J} e^{-2\pi \ii (j\cdot k) /n} z_j,\quad \forall k \in K,\\
z_j = (F^{-1} \hat{z})_j = \dfrac{1}{n^{d/2}}\sum_{k\in K} e^{+2\pi \ii (j\cdot k) /n} \hat{z}_k,\quad \forall j \in J.\\
\end{gather*}

The pseudospectral method discretizes the minus Laplacian operator
with
\begin{align*}
L := F^{-1} \tdiag(4\pi^2 |k|^2)_{k\in K} F
\end{align*}
which results in the discretized equation of \eqref{eqn:model}
\begin{equation}
\label{eqn:dmodel}
(L + \tdiag(v)) u = f,
\end{equation}
where, for example, $f = \{f_j = f(hj): j\in J\}$ and $v = \{v_j = v(hj): j\in J\}$
are the discrete arrays generated from sampling the values of $f(x)$ and $v(x)$ on the grid $hJ$,
respectively. $u = \{u_j : j\in J\}$ is the numerical solution on
$hJ$, where $u_j$ stands for an approximation of $u(hj)$. In what follows, we will use lower case letters to denote
discrete arrays on grids $J$ and $K$, which should not cause any
ambiguity.

\section{Brief Review of the Sparsifying Preconditioner}
\label{sec:review}
In this section we use a simplified version of the sparsifying
preconditioner in \cite{spspd} to review the main idea. While Equation \eqref{eqn:dmodel} gives an accurate dispersion relation, the numerical stencil is not local, which makes the sparse direct methods no longer applicable. The sparsifying preconditioner addresses this issue by approximating \eqref{eqn:dmodel} with a carefully designed sparse system, as we shall see in what follows.

Define $s:= (1/ |J|) \sum_{j\in J} v_j$ as the average of $v$ (a
scalar that is often quite negative in interesting cases, $s\sim
O(n^2)$) and $q := v - s$ as the zero-mean shift of $v$. Then
\eqref{eqn:dmodel} can be rewritten as
\begin{equation}
\label{eqn:dsmodel}
(L + s + \tdiag(q)) u = f.
\end{equation}
The reason why we want $q$ to have zero-mean will be explained
below.

We assume without loss of generality that $(L+s)$ is invertible,
otherwise we perturb $s$ by a small shift. The inverse of $(L+s)$, which is a highly indefinite matrix, is given by
\begin{gather*}
G := (L+s)^{-1} = F^{-1}\tdiag\left(\dfrac{1}{4\pi^2 |k|^2 + s}\right)_{k\in K} F,
\end{gather*}
which can be applied efficiently using FFT. Applying $G$ to both sides
of \eqref{eqn:dsmodel} gives
\begin{equation}
\label{eqn:G}
(I + G\, \tdiag(q)) u = Gf := g.
\end{equation}
The main motivation of the sparsifying preconditioner is that, $G$ is a Green's matrix induced from a partial differential equation and the operators in the equation are local. Suppose the step size $h$ is small enough, then if we discretize the equation with standard second order central difference scheme and we denote the resulting discrete operator as $A$, which is a sparse matrix, we will have that $AG \approx I$. Though the central difference scheme is not accurate enough when we only have a small number of points per oscillation, it is still reasonable to seek for some sparse matrix $Q$ as a ``sparse discretization'' of the equation, such that $QG$ is also sparse approximately, and by applying $Q$ on both sides of \eqref{eqn:G} we
get an approximately sparse system. The main task now is to find a $Q$
such that
\begin{enumerate}
\item
$Q$ is sparse and local,
\item
$QG$ is approximately sparse and local.
\end{enumerate}
Here, by ``sparse and local'', we mean that the non-zero elements in row $j$ only involve the nearby neighbors of $j$. If we could find such a $Q$ and denote $C \approx QG$ as the sparse
approximation, then applying $Q$ to \eqref{eqn:G} gives
\begin{equation}
\label{eqn:QG}
(Q+QG \, \tdiag(q)) u = Qg,
\end{equation}
which can be approximated by the sparse system
\begin{equation}
\label{eqn:QC}
(Q+C\, \tdiag(q)) \tilde{u} = Qg,
\end{equation}
where $\tilde{u}$ stands for an approximation of $u$. The sparse matrix $(Q+C\, \tdiag(q))$ can be inverted by the nested dissection algorithm which gives rise to an efficient preconditioner $(Q + C\, \tdiag(q))^{-1} Q$ to solve \eqref{eqn:G}.

To find such $Q$, we introduce the notation of neighborhood $\mu_j$ for
each $j\in J$
\begin{align*}
{\boldsymbol{\mu}}_j &:= \{i: \|i- j\|_\infty \le 1 \},
\end{align*}
which is the set containing $j$ and its nearest neighbors in
$l_\infty$ norm. Then the requirements for $Q$ can be formulated as
\begin{enumerate}
\item
$Q[j, {\boldsymbol{\mu}}_j^c] = 0$,
\item
$(QG)[j, {\boldsymbol{\mu}}_j^c] = Q[j, :] G[:, {\boldsymbol{\mu}}_j^c] = Q[j, {\boldsymbol{\mu}}_j] G[{\boldsymbol{\mu}}_j, {\boldsymbol{\mu}}_j^c] \approx 0$,
\end{enumerate}
where the notation $[\cdot,\cdot]$ means the submatrix of certain rows and columns. For example, $Q[j,{\boldsymbol{\mu}}_j^c]$ is the submatrix of $Q$ restricted to row $j$ and columns in ${\boldsymbol{\mu}}_j^c$. The superscript $c$ stands for complement and ${\boldsymbol{\mu}}_j^c = J \setminus {\boldsymbol{\mu}}_j$.

To find a proper choice for $Q[j,{\boldsymbol{\mu}}_j]$, we consider the following optimization
problem:
\begin{equation*}
\min_{\|\alpha\|_2 = 1} \| \alpha^T G[{\boldsymbol{\mu}}_j,{\boldsymbol{\mu}}_j^c]\|_2,
\end{equation*}
where $\alpha \in \bR^{|{\boldsymbol{\mu}}_j| \times 1}$ is a column vector. Let
$G[{\boldsymbol{\mu}}_j,{\boldsymbol{\mu}}_j^c]= U \Sigma V^T$ be the singular value decomposition. Then the optimal
solution for $\alpha$ is given by
\begin{align*}
\alpha &:= U[:,|{\boldsymbol{\mu}}_j|],
\end{align*}
where $|\cdot|$ means cardinality and $\alpha$ is the left singular vector corresponding to
the smallest singular value. We set
\begin{align*}
Q[j,{\boldsymbol{\mu}}_j]:= \alpha^T,
\end{align*}
then $\|(QG)[j,{\boldsymbol{\mu}}_j^c]\|_2$ is minimized and the optimal value is the smallest singular value. We expect the smallest singular value to be small because the partial differential equation itself implies that there should exist some local discretization to cancel the off diagonal elements of the Green's matrix approximately.

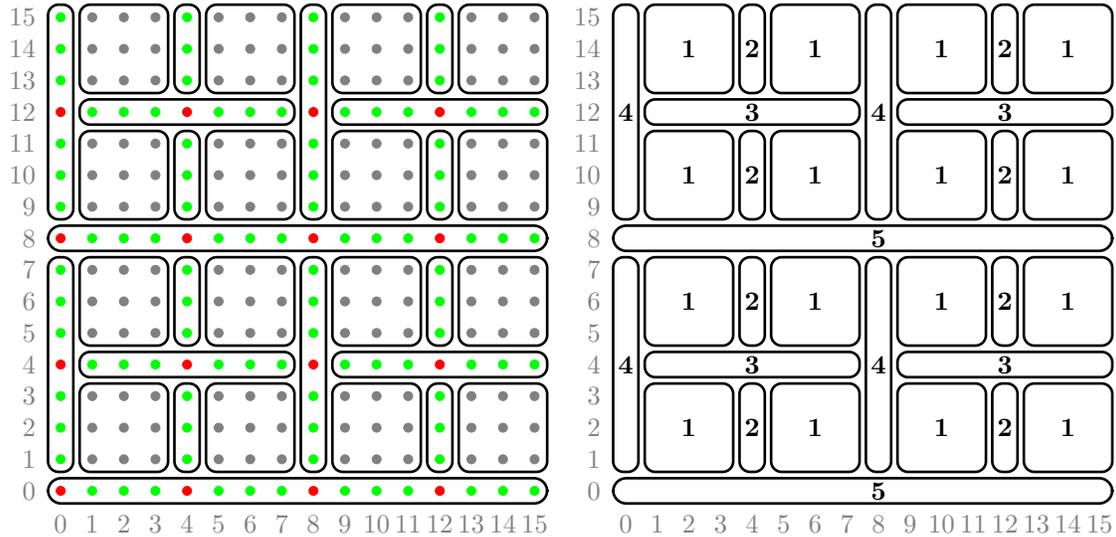
\begin{figure}[ht!]
\centering
\begin{tikzpicture}
[scale=0.42]
\foreach \x in {0,...,15}
\draw
[color = gray]
(\x,-0.5)node[below]{$\x$}
(-0.5,\x)node[left]{$\x$}
;

\foreach \x in {1,2,3,5,6,7,9,10,11,13,14,15}
\foreach \y in {1,2,3,5,6,7,9,10,11,13,14,15}
\filldraw[gray] (\x,\y) circle (4pt)
;
\foreach \x in {0,4,8,12}
\foreach \y in {1,2,3,5,6,7,9,10,11,13,14,15}
\filldraw[green] (\x,\y) circle (4pt)
;
\foreach \x in {1,2,3,5,6,7,9,10,11,13,14,15}
\foreach \y in {0,4,8,12}
\filldraw[green] (\x,\y) circle (4pt)
;
\foreach \x in {0,4,8,12}
\foreach \y in {0,4,8,12}
\filldraw[red] (\x,\y) circle (4pt)
;
\foreach \x in {-0.5}
\foreach \y in {-0.5,7.5}
\draw[rounded corners=5pt, line width = 1pt](\x+0.1,\y+0.1)rectangle(\x-0.1+16.0,\y-0.1+1.0)
;
\foreach \x in {-0.5,7.5}
\foreach \y in {0.5,8.5}
\draw[rounded corners=4pt, line width = 1pt](\x+0.1,\y+0.1)rectangle(\x-0.1+1.0,\y-0.1+7.0)
;
\foreach \x in {0.5,8.5}
\foreach \y in {3.5,11.5}
\draw[rounded corners=4pt, line width = 1pt](\x+0.1,\y+0.1)rectangle(\x-0.1+7.0,\y-0.1+1.0)
;
\foreach \x in {3.5,11.5}
\foreach \y in {0.5,4.5,8.5,12.5}
\draw[rounded corners=4pt, line width = 1pt](\x+0.1,\y+0.1)rectangle(\x-0.1+1.0,\y-0.1+3.0)
;
\foreach \x in {0.5,4.5,8.5,12.5}
\foreach \y in {0.5,4.5,8.5,12.5}
\draw[rounded corners=4pt, line width = 1pt](\x+0.1,\y+0.1)rectangle(\x-0.1+3.0,\y-0.1+3.0)
;
\end{tikzpicture}
\begin{tikzpicture}
[scale=0.42]
\foreach \x in {0,...,15}
\draw
[color = gray]
(\x,-0.5)node[below]{$\x$}
(-0.5,\x)node[left]{$\x$}
;
\foreach \x in {-0.5}
\foreach \y in {-0.5,7.5}
\draw[rounded corners=5pt, line width = 1pt](\x+0.1,\y+0.1)rectangle(\x-0.1+16.0,\y-0.1+1.0)
;
\foreach \x in {-0.5,7.5}
\foreach \y in {0.5,8.5}
\draw[rounded corners=4pt, line width = 1pt](\x+0.1,\y+0.1)rectangle(\x-0.1+1.0,\y-0.1+7.0)
;
\foreach \x in {0.5,8.5}
\foreach \y in {3.5,11.5}
\draw[rounded corners=4pt, line width = 1pt](\x+0.1,\y+0.1)rectangle(\x-0.1+7.0,\y-0.1+1.0)
;
\foreach \x in {3.5,11.5}
\foreach \y in {0.5,4.5,8.5,12.5}
\draw[rounded corners=4pt, line width = 1pt](\x+0.1,\y+0.1)rectangle(\x-0.1+1.0,\y-0.1+3.0)
;
\foreach \x in {0.5,4.5,8.5,12.5}
\foreach \y in {0.5,4.5,8.5,12.5}
\draw[rounded corners=4pt, line width = 1pt](\x+0.1,\y+0.1)rectangle(\x-0.1+3.0,\y-0.1+3.0)
;
\foreach \x in {2,6,10,14}
\foreach \y in {2,6,10,14}
\draw
(\x,\y)node{$\mathbf{1}$}
;
\foreach \x in {4,12}
\foreach \y in {2,6,10,14}
\draw
(\x,\y)node{$\mathbf{2}$}
;
\foreach \x in {4,12}
\foreach \y in {4,12}
\draw
(\x,\y)node{$\mathbf{3}$}
;
\foreach \x in {0,8}
\foreach \y in {4,12}
\draw
(\x,\y)node{$\mathbf{4}$}
;
\foreach \x in {8}
\foreach \y in {0,8}
\draw
(\x,\y)node{$\mathbf{5}$}
;
\end{tikzpicture}
\caption{This figure shows a $16\times 16$ example of the nested dissection algorithm. Left: The unknowns are grouped hierarchically. The gray points are the box-points in the nested dissection algorithm, the green ones are the edge-points and the red ones are the vertex-points. Right: the hierarchical elimination order is shown for each group. We first eliminate the box points $\mathbf{1}$ to their boundary neighbors and then eliminate the boundary points $\mathbf{2}$ to their remaining neighbors and so on so forth.}
\label{fig:nd}
\end{figure}

Once $Q$ is ready, we set $C$ as the truncation of $QG$ by
\begin{align*}
C[j,{\boldsymbol{\mu}}_j] &:= Q[j, {\boldsymbol{\mu}}_j] G[{\boldsymbol{\mu}}_j,{\boldsymbol{\mu}}_j],\\
C[j,{\boldsymbol{\mu}}_j^c] &:= 0.
\end{align*}
As a result, the matrix $P := Q + C\, \tdiag(q)$ has the same sparsity pattern as $Q$. By exploiting this sparsity, one can apply the nested dissection algorithm which reorders the unknowns hierarchically to minimize the elimination cost to solve the sparse system \eqref{eqn:QC}. Figure \ref{fig:nd} gives an example of the nested dissection algorithm in 2D.

We would like to point out that $\alpha, Q[j,{\boldsymbol{\mu}}_j]$ and
$C[j,{\boldsymbol{\mu}}_j]$ do not depend on $j$ due to the translational invariance
of $G$. Hence one only needs to perform the SVD and calculate the
stencils $Q[j,{\boldsymbol{\mu}}_j]$ and $C[j,{\boldsymbol{\mu}}_j]$ just once.

%

\section{Localized Sparsifying Preconditioner}
\label{sec:localized}
This section introduces the localized sparsifying preconditioner,
which is based on the sparsifying preconditioner in the previous
section but with two major improvements discussed below.

\subsection{Using the local potential information}
In Section \ref{sec:review}, the scalar shift $s$ is chosen to be the
average of $v$ such that $q$ has a zero-mean. The reason for choosing
$s$ to be the average is that, the error introduced in \eqref{eqn:QC}
comes from the truncation of $QG$. The truncated part, which should
have been multiplied by the elements in $q$, is neglected in
\eqref{eqn:QC}, and if $q$ is small in magnitude, then the error
introduced is expected to be small. The algorithm chooses $s$ such
that $q$ has a zero-mean to make the residual relatively small with a
single shift. However, it does not eliminate all the errors. When
$v(x)$ has large variations, the residual array $q$ could still be
large.

\vspace{0.1in}
\begin{figure}[ht!]
  \centering
  \begin{overpic}
    [width=0.9\textwidth]{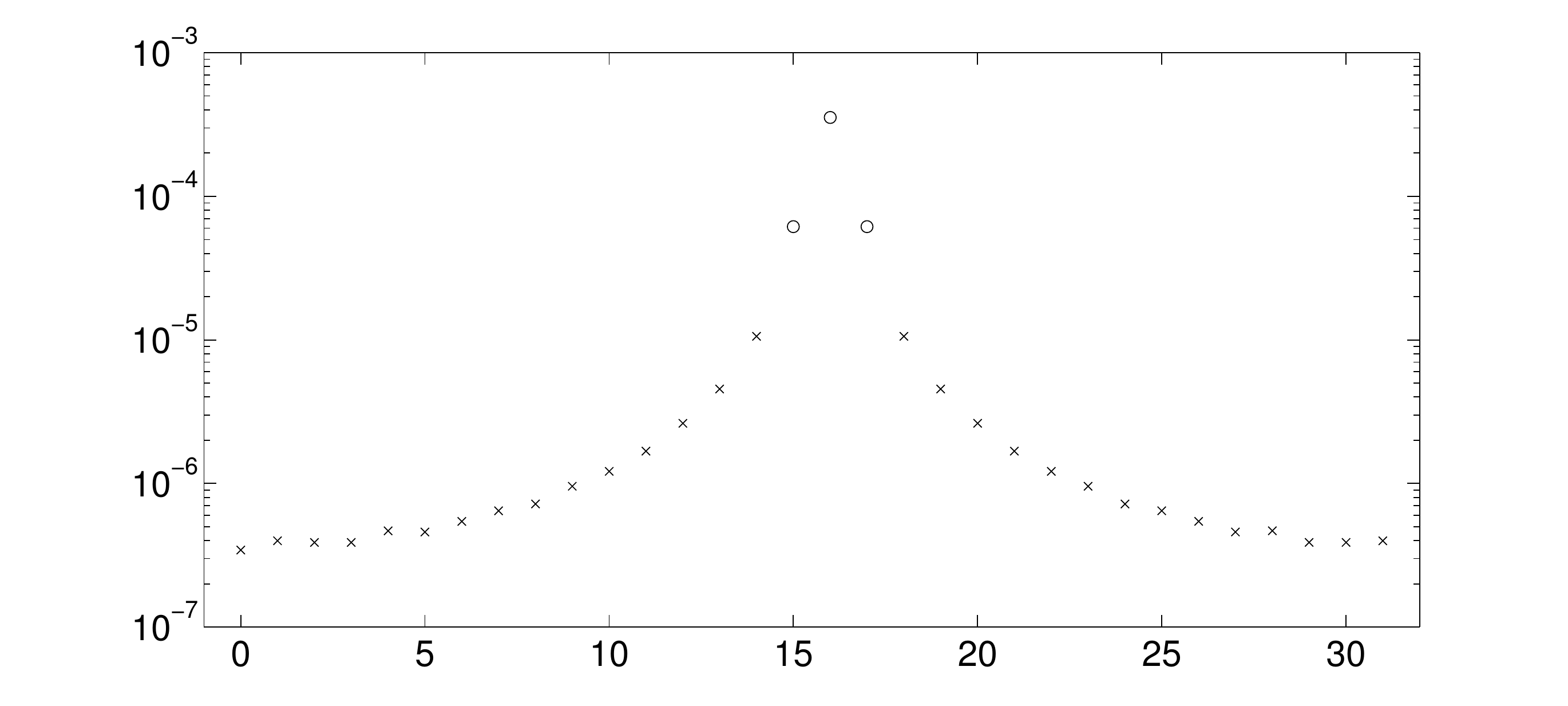}
    \put(3,45){$|QG[16,j_c]|$}
    \put(89,2){$j_c$}
  \end{overpic}
  \caption{This figure shows the magnitude of the $j$-th row of $QG$
    in logarithmic scale, where $n=32, j = 16$, and the shift $s = -62
    \pi^2$. The x-axis is the column index $j_c$ ranging from $0$ to
    $31$ and the y-axis is $|QG[16,j_c]|$. The entries with marks
    ``o'' are reserved and the ones with ``x'' are truncated. We can
    see the decay of the truncated entries as $j_c$ gets far away from
    the center position $j=16$.}
\label{fig:QGrow}
\end{figure}

Taking a closer look at the $j$-th row of $QG$, we find that the truncated
elements tend to be big near $j$. Figure \ref{fig:QGrow} shows a row
of $QG$ in the 1D case, which indicates a decaying trend of the
elements of the $j$-th row of $QG$ when the column index gets far away
from $j$. It tells us that, if we choose the shift $s$ such that $q$
is small near $j$, then the truncated elements in the $j$-th row of
$QG\, \tdiag(q)$ will be more likely to be diminished, since the dominating
elements in $QG$ are multiplied with small elements in $q$. This
suggests that, it is helpful to set the shift $s$ to be close to $v_j$
to reduce the error introduced by truncating the $j$-th row of
$QG\, \tdiag(q)$. Ideally, it would be nice if, for each $j$, one could use the
Green's matrix $G_{v_j} :=(L + v_j)^{-1}$ to compute the SVD and set
the stencils by
\begin{align*}
G_{v_j}[{\boldsymbol{\mu}}_j,{\boldsymbol{\mu}}_j^c] &= U \Sigma V^T,\\
\alpha_{v_j} &:= U[:,|{\boldsymbol{\mu}}_j|],\\
Q[j, {\boldsymbol{\mu}}_j] &:= \alpha_{v_j}^T,\\
C[j,{\boldsymbol{\mu}}_j] &:= Q[j, {\boldsymbol{\mu}}_j] G_{v_j}[{\boldsymbol{\mu}}_j,{\boldsymbol{\mu}}_j],
\end{align*}
in which case the truncated part tends to be much smaller.

However, computing each row of $Q$ and $C$ with a unique shift is expensive. To
save computational cost, a list of possible shifts is created in
advance. Then our method only computes the local stencils
corresponding to the shifts in the list and assigns the stencil at
position $j$ to the shift closest to $v_j$. More specifically, we
first choose a set $S$ which contains a list of shifts $s$ distributed
evenly in the range of $v$. The method of choosing $S$ will be
discussed later. Then for each $s\in S$, one computes $\alpha_s$ by
\begin{align*}
G_{s} &:= (L + s)^{-1},\\
G_{s}[{\boldsymbol{\mu}}_j,{\boldsymbol{\mu}}_j^c] &= U\Sigma V^T,\\
\alpha_{s} &:= U[:,|{\boldsymbol{\mu}}_j|].
\end{align*}
Notice that the value of $\alpha_s$ does not depend on $j$ due to the
translational invariant property of $G_s$.

After calculating $\alpha_s$ for each $s$, our method sets a shift $s_j$ for
each $j$ to be the shift $s$ closest to $v_j$ in $S$
\begin{align*}
s_j := \min \{|s-v_j|: s\in S\},
\end{align*}
which means $s_j$ serves as an approximation of the local shift $v_j$ hence it is location dependent. However, the the singular vector $\alpha_{s_j}$ only needs to be computed once for different locations sharing the same shift approximation. Thus the number of Green's matrices and SVDs that need to be formed only depends on the range of the potential shift $v$, not on the number of discrete points. That saves us the computational cost, especially in 3D case.

With $\alpha_{s_j}$, one computes the following stencils for each $j$
\begin{align*}
Q[j, {\boldsymbol{\mu}}_j] &:= \alpha_{s_j}^T,\\
C[j,{\boldsymbol{\mu}}_j] &:= Q[j, {\boldsymbol{\mu}}_j] G_{s_j}[{\boldsymbol{\mu}}_j,{\boldsymbol{\mu}}_j],\\
P[j, {\boldsymbol{\mu}}_j] &:= Q[j, {\boldsymbol{\mu}}_j] + C[j, {\boldsymbol{\mu}}_j] \, \tdiag(v[{\boldsymbol{\mu}}_j] - s_j)
\end{align*}
where $v[{\boldsymbol{\mu}}_j]$ is the array $v$ restricted to ${\boldsymbol{\mu}}_j$.  Now
multiplying $G_{s_j}$ to \eqref{eqn:dmodel} on both sides gives
\begin{align}
\label{eqn:GNj}
(I + G_{s_j}\, \tdiag(v-s_j)) u = G_{s_j} f.
\end{align}
Next multiplying by $Q[j, {\boldsymbol{\mu}}_j]$ to the rows indexed by ${\boldsymbol{\mu}}_j$ in
\eqref{eqn:GNj} and truncating the elements not in ${\boldsymbol{\mu}}_j$ gives rise
to
\begin{align*}
P[j, {\boldsymbol{\mu}}_j] u[{\boldsymbol{\mu}}_j] \approx C[j, {\boldsymbol{\mu}}_j] f[{\boldsymbol{\mu}}_j].
\end{align*}
Assembling all the approximating equations for different positions $j$ results in
the following equation
\begin{equation}
\label{eqn:PC}
P \tilde{u} = C f,
\end{equation}
where $\tilde{u}$ serves as an approximation of $u$. By applying the nested dissection algorithm, one can use $P^{-1} C$ as an efficient preconditioner for \eqref{eqn:dmodel}. We see that the local potential information is taken into consideration to build the
preconditioner, which is the main advantage of this new approach.


Now let us discuss how to choose $S$. Denote $v_{\min}$ and $v_{\max}$
as the minimum and maximum value of $v$ respectively. To approximate
all local potential values in the interval $[v_{\min}, v_{\max}]$ by a
minimal distance, our method distributes the shifts in $S$ evenly in
$[v_{\min}, v_{\max}]$. Another important point for choosing the shift
$s$ is to avoid resonance, which means, we do not want $(L+ s)$ to be
singular. A simple way to avoid resonance is to choose $s$ to have the
form $(4m + 2)\pi^2$ where $m \in \bZ$, since the eigenvalues of $L$
are all multiples of $4 \pi^2$. For the size of the shift list $S$,
the numerical tests in Section \ref{sec:numerical} show that it
suffices to set $|S| = O(n)$ to get an ideal iteration number.

\subsection{Computing local stencils using FFT}
As pointed out in \cite{spspd}, building compact stencils that only
involve the nearest neighbors fails to give a preconditioner accurate
enough. \cite{spspd} solves this issue by treating the points in the
same leaf box of the nested dissection algorithm as a whole and
increasing the size of the leaf box as the problems size grows. This
paper adopts a different approach.

For each $s \in S$, we not only build the stencil involving the
nearest neighbors but also the ones involving neighbors at a larger
distance, in order to get more accurate local schemes. Then, for each
position $j$, our method picks the most accurate stencil that involves
the largest possible neighborhood such that the nested dissection
algorithm can still be applied. More specifically, define
\[
{\boldsymbol{\mu}}_j^t := \{i : \|i - j\|_\infty \le t \}
\]
where $t$ controls the size of the neighborhood and $t$ is typically
set to be bounded by $4$. For example, $t = 1$ corresponds to the
compact stencil we discussed above. For each $(s,t)$ pair, we compute
\begin{align*}
  G_s[{\boldsymbol{\mu}}_j^t, ({\boldsymbol{\mu}}_j^t)^c] &= U \Sigma V^T,\\
  \alpha_s^t &:= U[:,|{\boldsymbol{\mu}}_j^t|],
\end{align*}
where $\alpha_s^t$ does not depend on $j$ due to the translational invariance of the Green's matrix, and we can simply use $G_s[{\boldsymbol{\mu}}_0^t, ({\boldsymbol{\mu}}_0^t)^c]$ to compute the SVD.

\begin{figure}[ht!]
\centering
\begin{tikzpicture}
[scale=0.6]
\draw
(-0.5,-0.5)grid(15.5,15.5)
;
\foreach \x in {0,...,15}
\draw
[color = gray]
(\x,-0.5)node[below]{$\x$}
(-0.5,\x)node[left]{$\x$}
;
\foreach \x in {0,...,15}
\foreach \y in {0,...,15}
\draw
(\x,\y)node[fill=gray]{$1$}
;
\foreach \x in {2,6,10,14}
\foreach \y in {2,6,10,14}
\draw
(\x,\y)node[fill=gray]{$2$}
;
\foreach \x in {0,...,15}
\foreach \y in {0,4,8,12}
\draw
(\x,\y)node[fill=green]{$1$}
;
\foreach \x in {0,4,8,12}
\foreach \y in {0,...,15}
\draw
(\x,\y)node[fill=green]{$1$}
;
\foreach \x in {2,6,10,14}
\foreach \y in {0,4,8,12}
\draw
(\x,\y)node[fill=green]{$2$}
;
\foreach \x in {0,4,8,12}
\foreach \y in {2,6,10,14}
\draw
(\x,\y)node[fill=green]{$2$}
;
\foreach \x in {0,4,8,12}
\foreach \y in {0,4,8,12}
\draw
(\x,\y)node[fill=red]{$4$}
;
\end{tikzpicture}

\caption{This figure shows a $16\times 16$ example of the choice of $t$ in 2D. The gray points are the box-points in the nested dissection algorithm, the green ones are the edge-points and the red ones are the vertex-points. The maximum choice of $t$ is marked out for each point. To satisfy the sparsity pattern requirement for the nested dissection algorithm, the interaction between two points cannot cross edges. For example, two box-points in different boxes cannot interact with each other. We see that for the box-points, the ones nearer to the box center have larger stencil supports. We also note that the vertex-points can have larger supports without breaking the sparsity requirement.}
\label{fig:tchoice}
\end{figure}
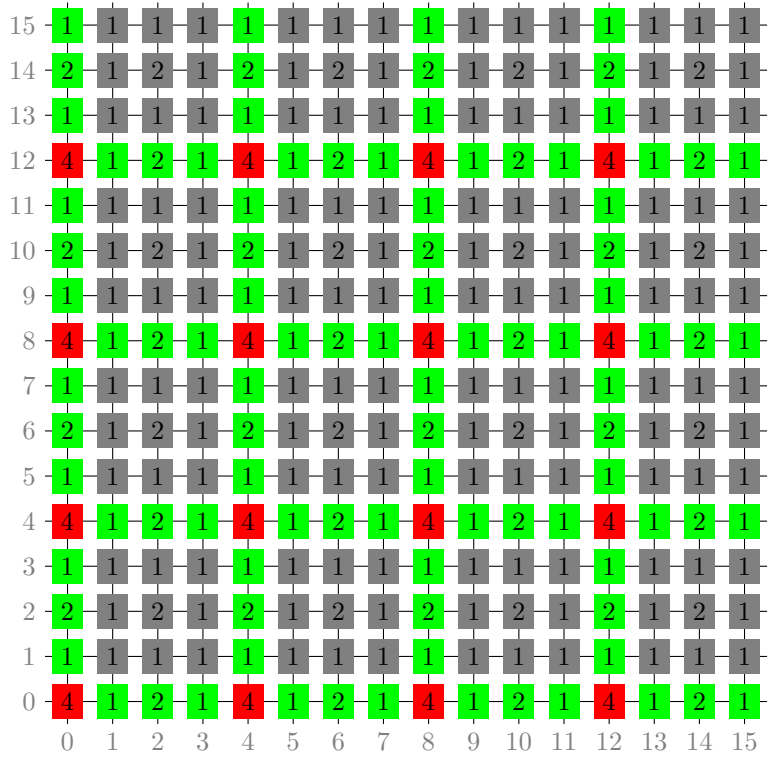

Now for each $j$, set $s = s_j$ and choose $t$ as big as possible such
that the sparsity pattern requirement for the nested dissection
algorithm is satisfied. Figure \ref{fig:tchoice} shows an example for
the choice of $t$ for different points, where a point nearer to the
center of a box has a larger $t$ value. In what follows, we denote
$t_j$ as this neighbor size chosen for the position $j$. Then the
stencils are given by
\begin{align*}
  Q[j,{\boldsymbol{\mu}}_j^{t_j}] &:= (\alpha_{s_j}^{t_j})^T,\\
  C[j,{\boldsymbol{\mu}}_j^{t_j}] &:= Q[j, {\boldsymbol{\mu}}_j^{t_j}] G_{s_j}[{\boldsymbol{\mu}}_j^{t_j},{\boldsymbol{\mu}}_j^{t_j}],\\
  P[j,{\boldsymbol{\mu}}_j^{t_j}] &:= Q[j, {\boldsymbol{\mu}}_j^{t_j}] + C[j, {\boldsymbol{\mu}}_j^{t_j}] \, \tdiag(v[{\boldsymbol{\mu}}_j^{t_j}] - s_j),
\end{align*}
where $v[{\boldsymbol{\mu}}_j^{t_j}]$ stands for restricting the vector $v$ to the
index set ${\boldsymbol{\mu}}_j^{t_j}$.

Let us discuss the computational cost of building the stencils. First,
the cost of forming the Green's matrix $G_s$ is $O(N \log N)$ since we
only need to compute a single column of $G_s$ via FFT. Next, if we use
the SVD of $G_s[{\boldsymbol{\mu}}_0^t,({\boldsymbol{\mu}}_0^t)^c]$ directly to compute $\alpha_{s}^{t}$, the
cost will be $O(t^{2d} N)$, which can dominate the $O(N \log N)$
complexity in practice even when $t$ is as small as $2$ or $3$. To
reduce the cost, consider the following identity
\begin{align}
G_s[{\boldsymbol{\mu}}_0^t,({\boldsymbol{\mu}}_0^t)^c] G_s[{\boldsymbol{\mu}}_0^t,({\boldsymbol{\mu}}_0^t)^c]^T
&= G_s[{\boldsymbol{\mu}}_0^t, :]G_s[{\boldsymbol{\mu}}_0^t,:]^T - G_s[{\boldsymbol{\mu}}_0^t, {\boldsymbol{\mu}}_0^t] G_s[{\boldsymbol{\mu}}_0^t, {\boldsymbol{\mu}}_0^t]^T\\
&:= A_1 - A_2.
\label{eqn:A}
\end{align}
In order to compute $\alpha_{s}^{t}$, one only needs to calculate the
eigenvalue decomposition of $(A_1 - A_2)$ and extract the eigenvector
corresponding to the smallest eigenvalue. A nice property of
\eqref{eqn:A} is that, the elements in $A_1$ are dot products of the
rows of $G_s$, which are essentially the values of the convolution of
the kernel vector with itself. To explain in details, denote $g_s$ as
the kernel vector of $G_s$ with periodic extension over $Z^d$. Then
\begin{align*}
G_s[j_1, j_2] = g_s[j_1 - j_2],\quad & \forall j_1,j_2 \in J,\\
g_s[j] = g_s[-j],\quad &\forall j \in J.
\end{align*}

The entries in $A_1$ have the form $G_s[j_1,:] G_s[j_2, :]^T$, which can be written as
\begin{align*}
G_s[j_1,:] G_s[j_2,:]^T &= \sum_{j \in J} g_s[j_1 - j] g_s [j - j_2]\\
&= \sum_{j \in J} g_s[j_1 - j_2 - j] g_s [j]\\
&= (g_s \ast g_s) [j_1 - j_2].
\end{align*}

The convolution $g_s \ast g_s$ can be computed by FFT in $O(N \log N)$
steps. After computing the convolution, retrieving the entries in
$A_1$ costs only $O(t^{2d})$ steps. Retrieving the other entries in
\eqref{eqn:A} and computing the eigenvalue decomposition take
$O(t^{3d})$ steps so the total cost is $O(N \log N + t^{3d})$ for
computing the stencil $\alpha_{s}^{t}$ for each $(s, t)$ pair. Here $O(t^{3d})$
is negligible compared to $O(N \log N)$ in practical cases. In addition,
the computation of the convolution of $g_s$ with itself only needs to
be performed once for different $t$ values, so there is only little extra cost for
building stencils involving larger range of neighbor points.  The
overall computational costs are listed in Table \ref{tab:cost}.
\begin{table}[ht!]
\centering
\begin{tabular}{c|c|c|c}
\hline
 & $T_{\text{stencil}}$ & $T_{\text{NDsetup}}$ & $T_{\text{NDapp}}$ \\
 \hline
2D & $O(|S| (N \log N + t^{6}))$ & $O(N^{3/2} + b^4 N)$ & $O(N \log N + b^2 N)$\\
\hline
3D & $O(|S| (N \log N + t^{9}))$ & $O(N^2 + b^6 N)$ & $O(N^{4/3} + b^3 N)$\\
\hline
\end{tabular}
\caption{The time costs of the algorithm. $T_{\text{stencil}}$ is the cost of computing the stencils used for building $Q, C$ and $P$. $T_{\text{NDsetup}}$ is the setup cost of the nested dissection algorithm and $T_{\text{NDapp}}$ is the application cost per iteration. $N = n^d$ is the degree of freedom. $b$ is the size of the leaf box in the nested dissection algorithm. $t$ is the size of the largest neighborhood and $|S|$ is the size of the shift list.}
\label{tab:cost}
\end{table}

\section{Numerical Results}
\label{sec:numerical}
This section presents the numerical results for the Helmholtz equation
and the Schr\"odinger equation in 2D and 3D. The algorithm is
implemented in MATLAB and the tests are performed on a server with
four Intel Xeon E5-4640 CPUs and a max usage of 384 GB memory. The
preconditioner is combined with standard GMRES solver with relative
tolerance $10^{-6}$ and restart value $40$.  The notations in the
numerical tests are listed as follows.
\begin{itemize}
\item $\omega$ is the angular frequency.
\item $N=n^d$ is the number of unknowns.
\item $|S|$ is the size of the shift list.
\item $T_{\text{stencil}}$ is the time cost of computing the stencils in seconds.
\item $T_{\text{NDsetup}}$ is the setup cost of the nested dissection algorithm in seconds.
\item $N_{\text{iter}}$ is the iteration number.
\item $T_{\text{NDsolve}}$ is the solve cost of the nested dissection algorithm in seconds.
\end{itemize}

Our method computes the stencils for $t=1,2$ and sets $|S| = O(n)$. The
leaf box in the nested dissection algorithm $b$ is fixed to be
$8$. The right-hand side for each test is a Gaussian point source at the center of the domain.

\paragraph{Helmholtz Equation.}
For the Helmholtz equation, $v(x) = -(\omega/c(x))^2$ where $\omega$
is the angular frequency and $c(x)$ is the velocity field. Two
velocity fields in 2D are tested:
\begin{enumerate}  [(i).]
\item  \label{2DHi}  A constant background with a Gaussian profile in the center of the square.
\item  \label{2DHii}  A constant background with a cross shape profile.
\end{enumerate}
The results are given in Tables \ref{tab:2DHi} and \ref{tab:2DHii}.

\begin{table}[ht!]
  \centering
  \begin{tabular}{ccc|cc|cc}
    \hline
    $\omega/(2\pi)$ & $N$ & $|S|$ & $T_{\text{stencil}}$ & $T_{\text{NDsetup}}$ & $N_{\text{iter}}$ & $T_{\text{NDsolve}}$ \\
    \hline
$16$ & $64^2$ & $4$ & $1.89e-02$ & $5.25e-02$ & $6$ & $4.00e-02$ \\ 
$32$ & $128^2$ & $8$ & $2.71e-02$ & $2.27e-01$ & $5$ & $1.39e-01$ \\ 
$64$ & $256^2$ & $16$ & $8.30e-02$ & $8.91e-01$ & $6$ & $6.74e-01$ \\ 
$128$ & $512^2$ & $32$ & $3.87e-01$ & $4.23e+00$ & $6$ & $2.56e+00$ \\ 
    \hline
  \end{tabular}
  \includegraphics[width=0.425\textwidth]{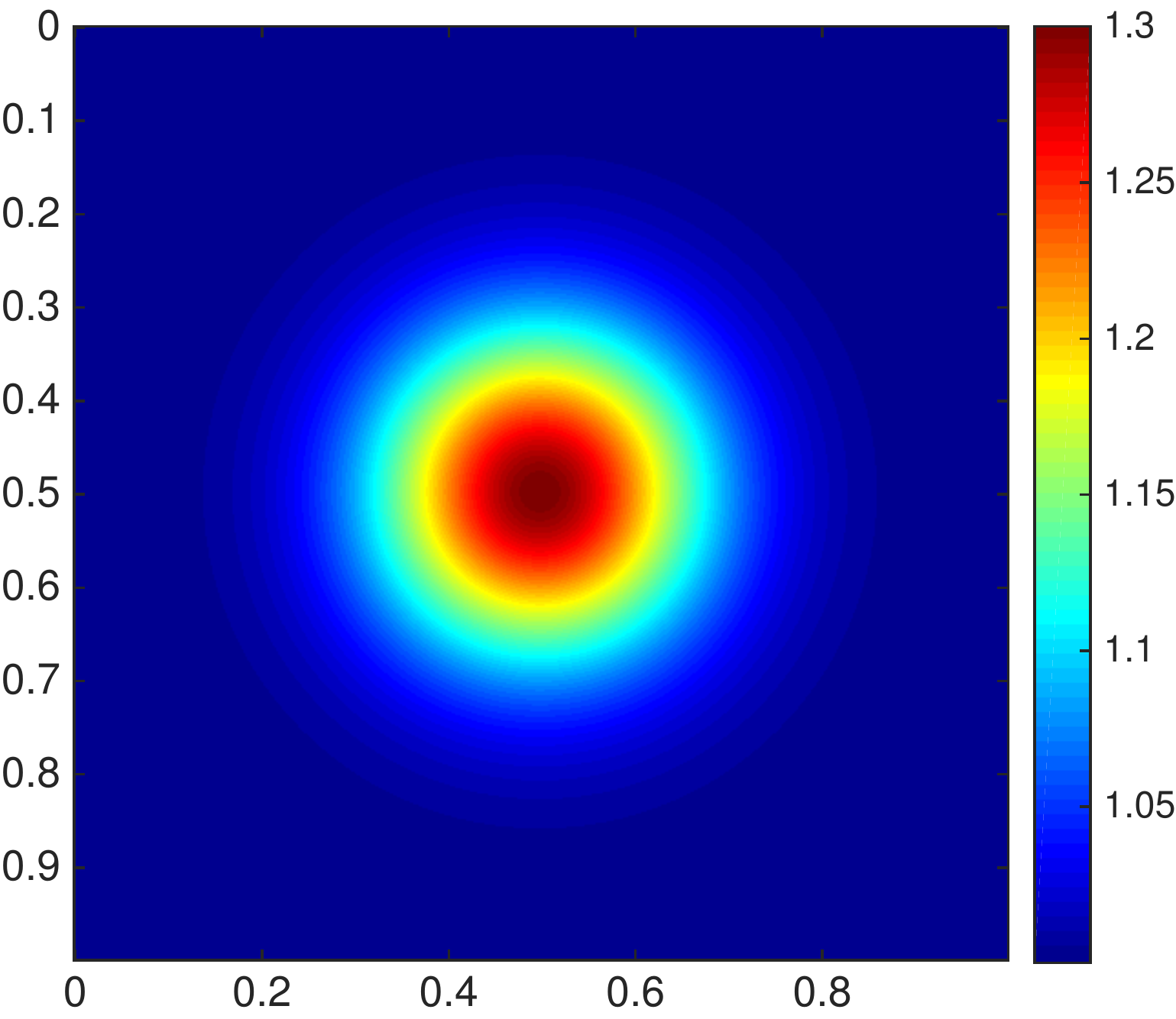}
  \includegraphics[width=0.425\textwidth]{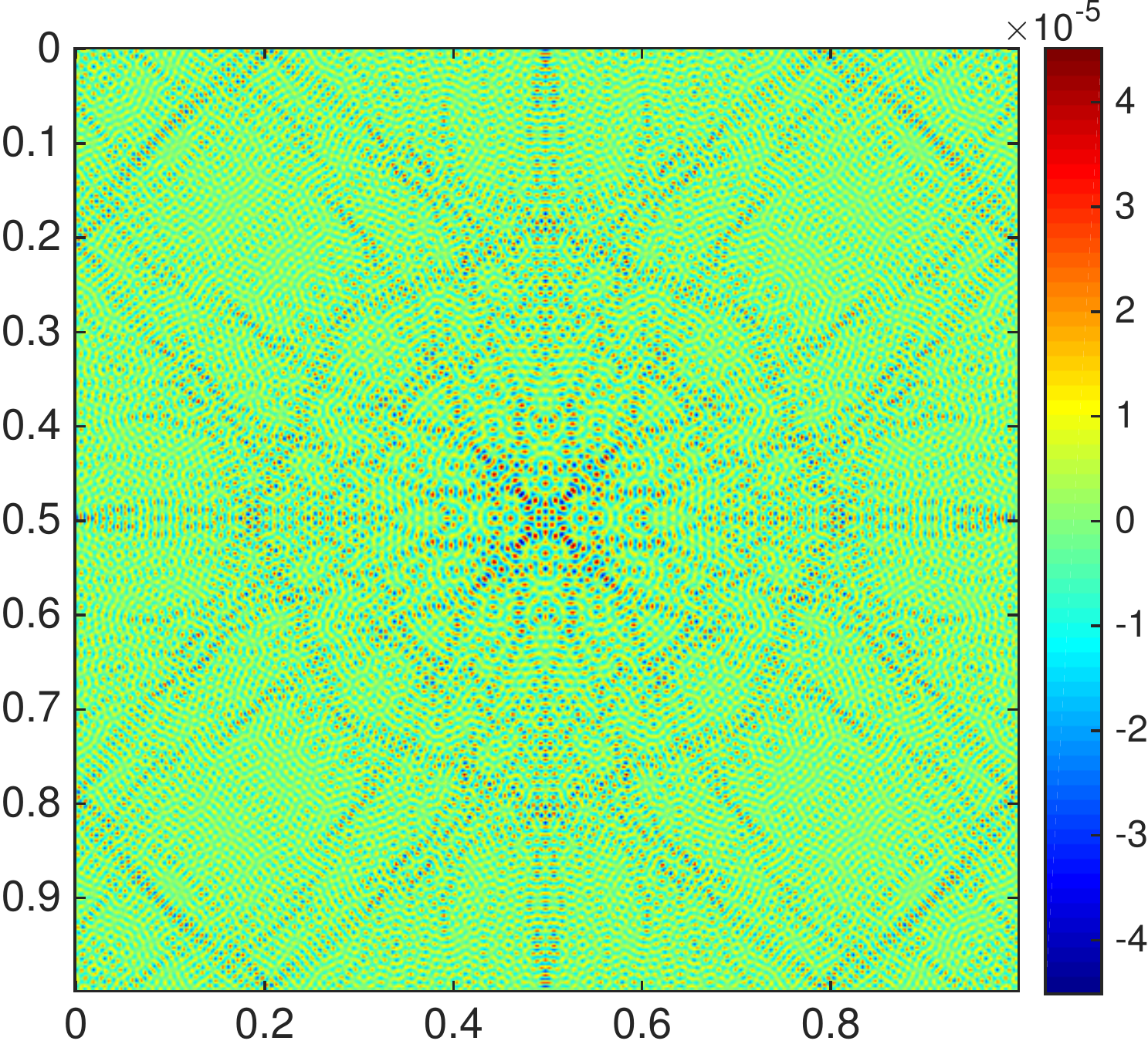}
  \caption{Results for velocity field (\ref{2DHi}) of the 2D Helmholtz equation. Top: numerical results. Bottom: $c(x)$ (left) and $u(x)$ (right) for the largest problem size.}
  \label{tab:2DHi}
\end{table}

\begin{table}[ht!]
  \centering
  \begin{tabular}{ccc|cc|cc}
    \hline
    $\omega/(2\pi)$ & $N$ & $|S|$ & $T_{\text{stencil}}$ & $T_{\text{NDsetup}}$ & $N_{\text{iter}}$ & $T_{\text{NDsolve}}$ \\
    \hline
$16$ & $64^2$ & $4$ & $8.39e-03$ & $5.23e-02$ & $6$ & $3.66e-02$ \\ 
$32$ & $128^2$ & $8$ & $3.89e-02$ & $2.12e-01$ & $6$ & $1.66e-01$ \\ 
$64$ & $256^2$ & $16$ & $7.66e-02$ & $9.10e-01$ & $6$ & $6.67e-01$ \\ 
$128$ & $512^2$ & $32$ & $3.62e-01$ & $4.49e+00$ & $5$ & $2.51e+00$ \\ 
    \hline
  \end{tabular}
  \includegraphics[width=0.425\textwidth]{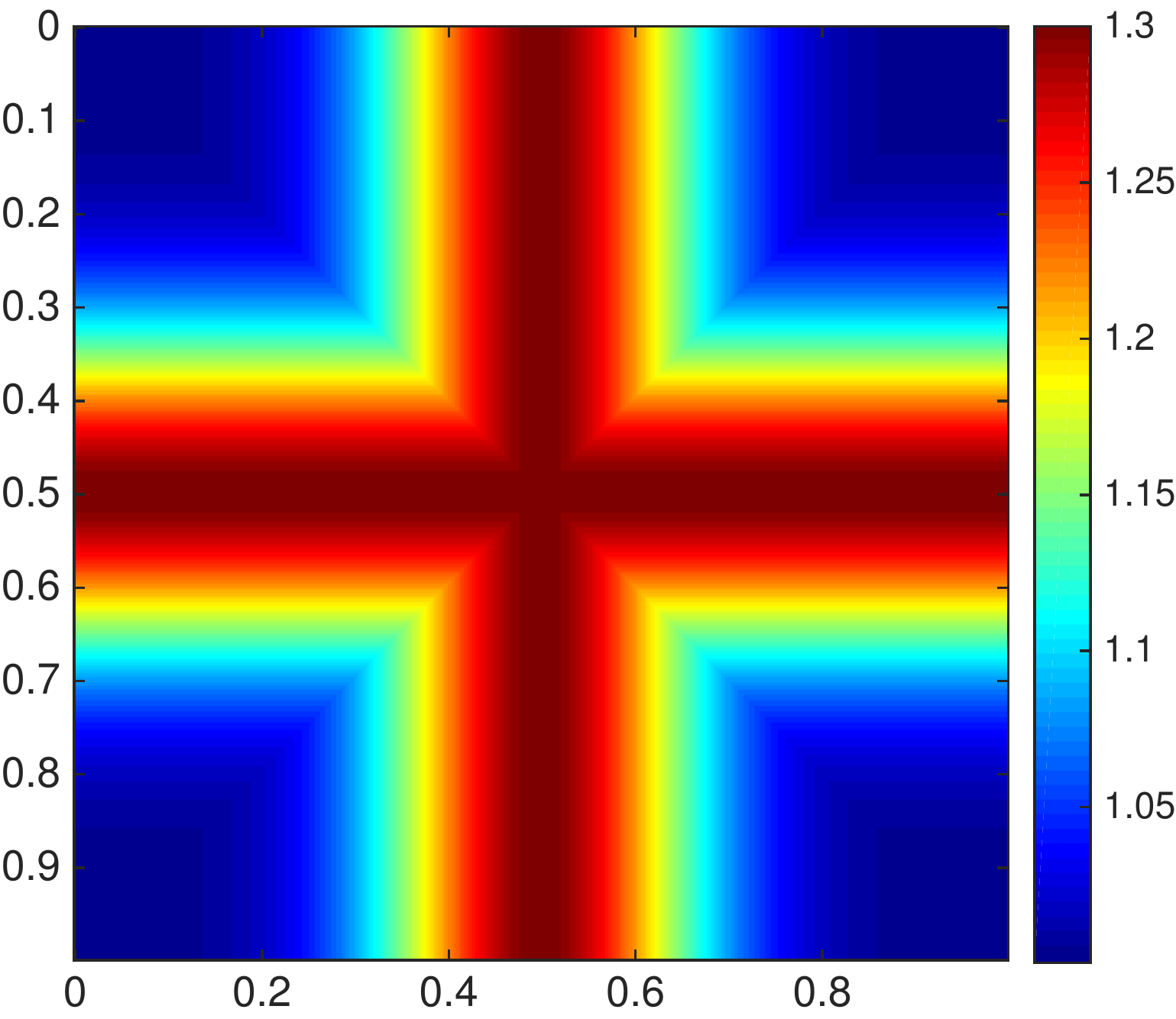}
  \includegraphics[width=0.425\textwidth]{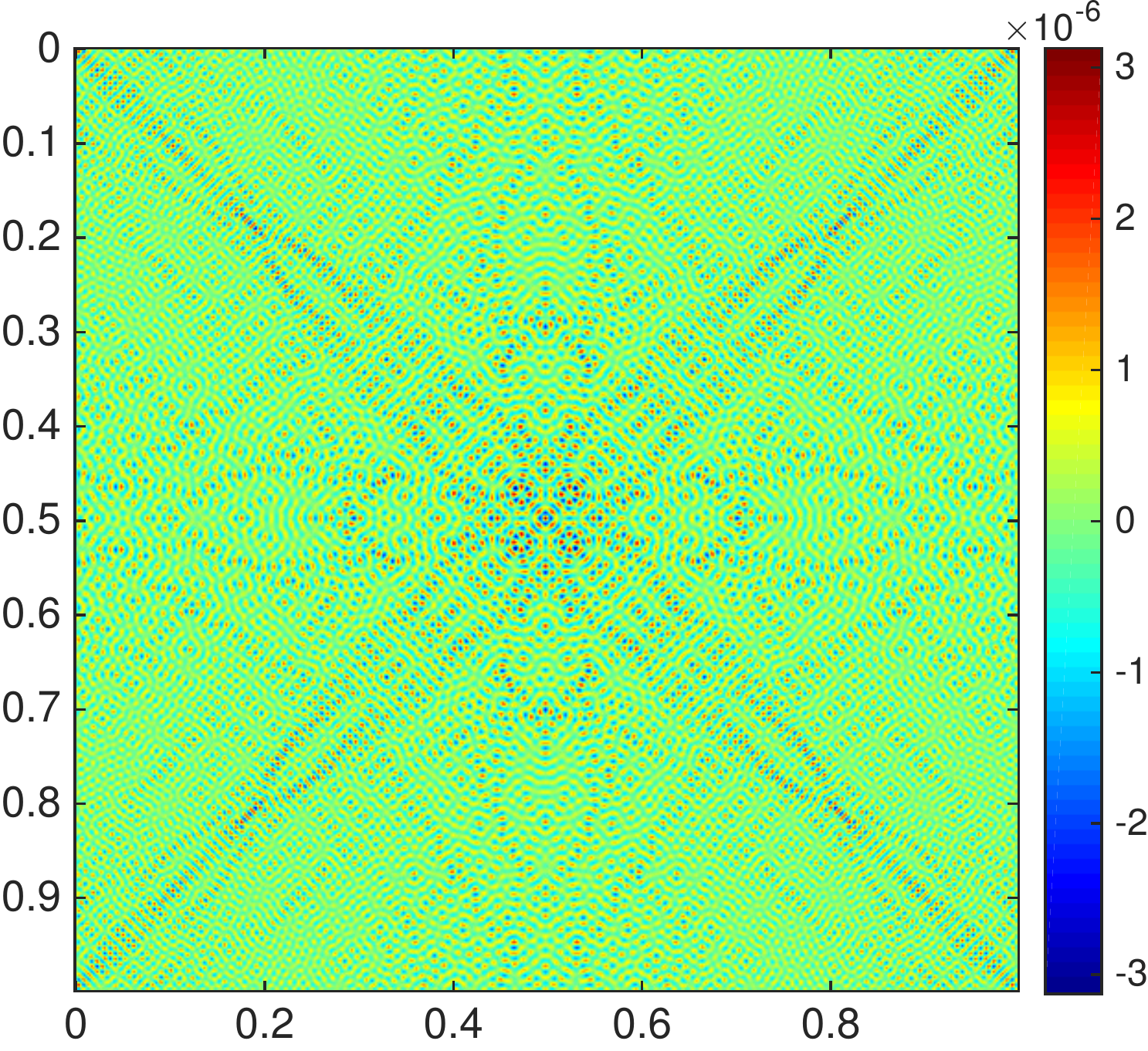}
  \caption{Results for velocity field (\ref{2DHii}) of the 2D Helmholtz equation. Top: numerical results. Bottom: $c(x)$ (left) and $u(x)$ (right) for the largest problem size.}
  \label{tab:2DHii}
\end{table}

Two similar tests are performed in 3D, where the velocity fields are
\begin{enumerate}[(i).]
\item \label{3DHi} A constant background with a Gaussian profile in
  the center of the cube.
\item \label{3DHii} A constant background with a cross shape profile,
  where three cylinders cross the domain along the three axes
  correspondingly.
\end{enumerate}
The results are summarized in Tables \ref{tab:3DHi} and \ref{tab:3DHii}.

\begin{table}[ht!]
  \centering
  \begin{tabular}{ccc|cc|cc}
    \hline
    $\omega/(2\pi)$ & $N$ & $|S|$ & $T_{\text{stencil}}$ & $T_{\text{NDsetup}}$ & $N_{\text{iter}}$ & $T_{\text{NDsolve}}$ \\
    \hline
$4$ & $16^3$ & $4$ & $4.80e-02$ & $2.55e-01$ & $6$ & $8.55e-02$ \\ 
$8$ & $32^3$ & $8$ & $8.37e-02$ & $6.20e+00$ & $6$ & $9.10e-01$ \\ 
$16$ & $64^3$ & $16$ & $3.73e-01$ & $1.70e+02$ & $6$ & $1.20e+01$ \\ 
$32$ & $128^3$ & $32$ & $3.53e+00$ & $8.72e+03$ & $9$ & $2.07e+02$ \\ 
    \hline
  \end{tabular}
  \includegraphics[width=0.425\textwidth]{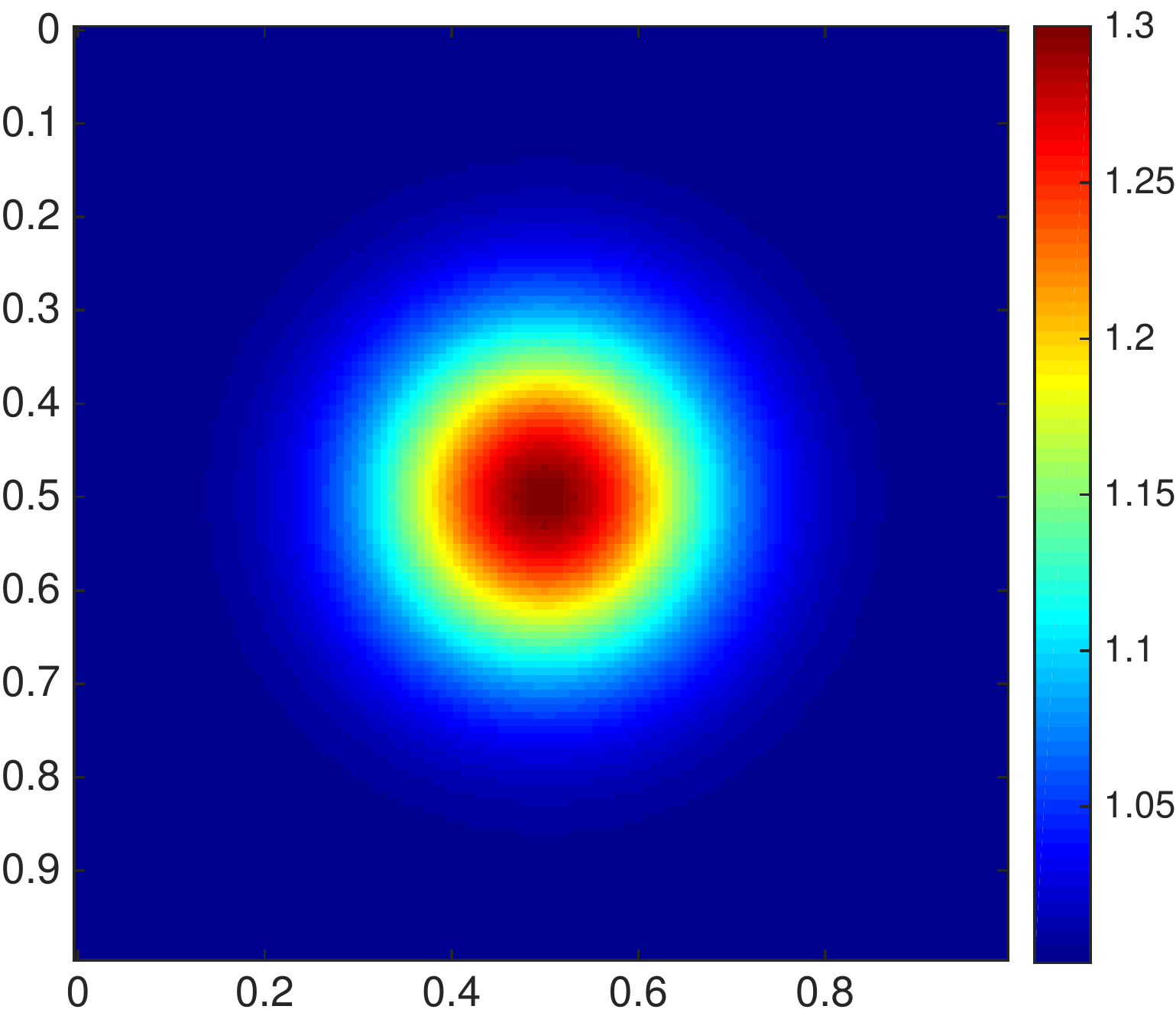}
  \includegraphics[width=0.425\textwidth]{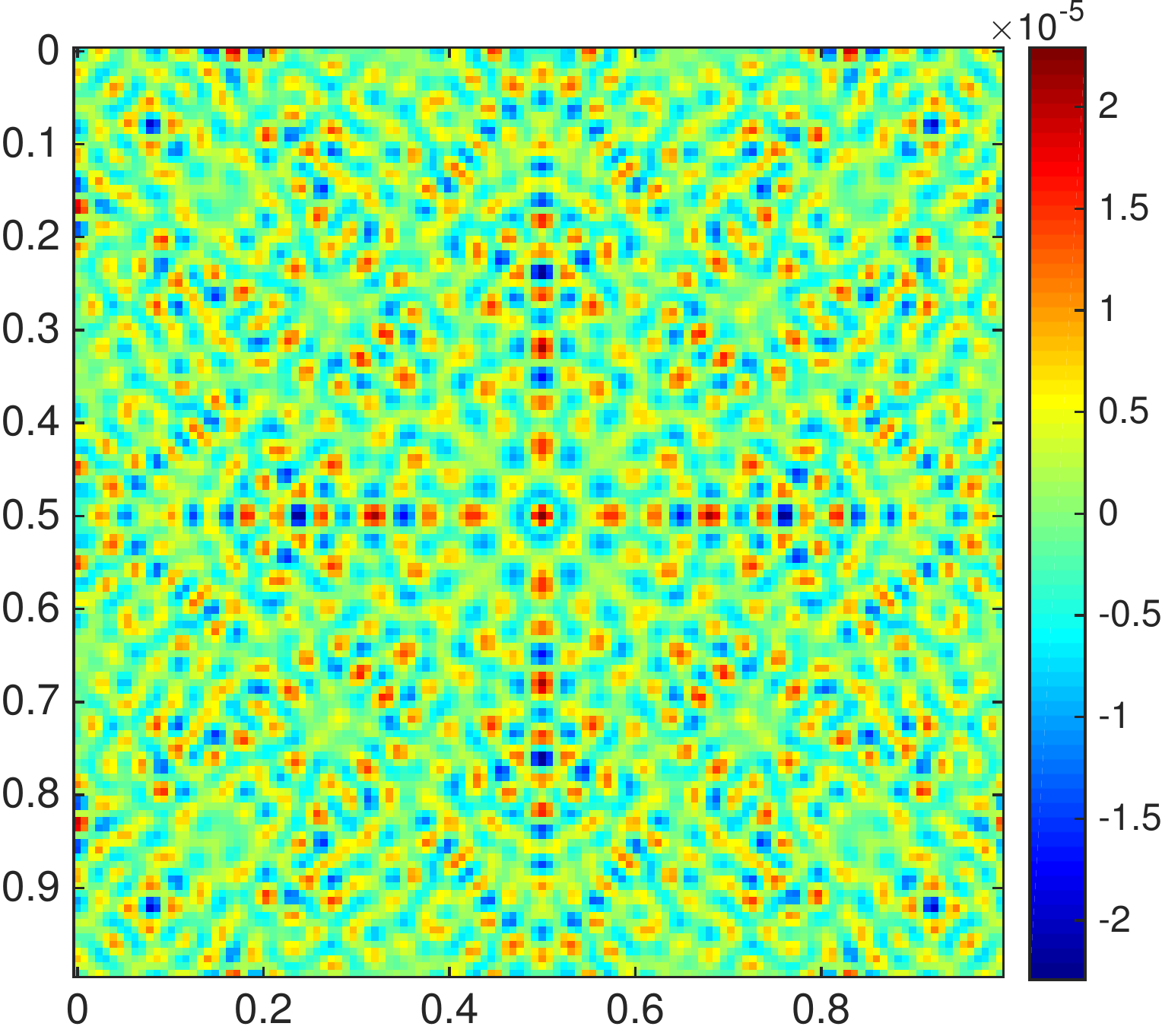}
  \caption{Results for velocity field (\ref{3DHi}) of the 3D Helmholtz equation. Top: numerical results. Bottom: $c(x)$ (left) and $u(x)$ (right) at $x_3 = 0.5$ for the largest problem size.}
  \label{tab:3DHi}
\end{table}

\begin{table}[ht!]
  \centering
  \begin{tabular}{ccc|cc|cc}
    \hline
    $\omega/(2\pi)$ & $N$ & $|S|$ & $T_{\text{stencil}}$ & $T_{\text{NDsetup}}$ & $N_{\text{iter}}$ & $T_{\text{NDsolve}}$ \\
    \hline
$4$ & $16^3$ & $4$ & $3.76e-02$ & $2.91e-01$ & $7$ & $9.51e-02$ \\ 
$8$ & $32^3$ & $8$ & $8.60e-02$ & $6.00e+00$ & $6$ & $9.93e-01$ \\ 
$16$ & $64^3$ & $16$ & $3.85e-01$ & $1.64e+02$ & $6$ & $1.12e+01$ \\ 
$32$ & $128^3$ & $32$ & $4.07e+00$ & $8.64e+03$ & $7$ & $1.47e+02$ \\ 
    \hline
  \end{tabular}
  \includegraphics[width=0.425\textwidth]{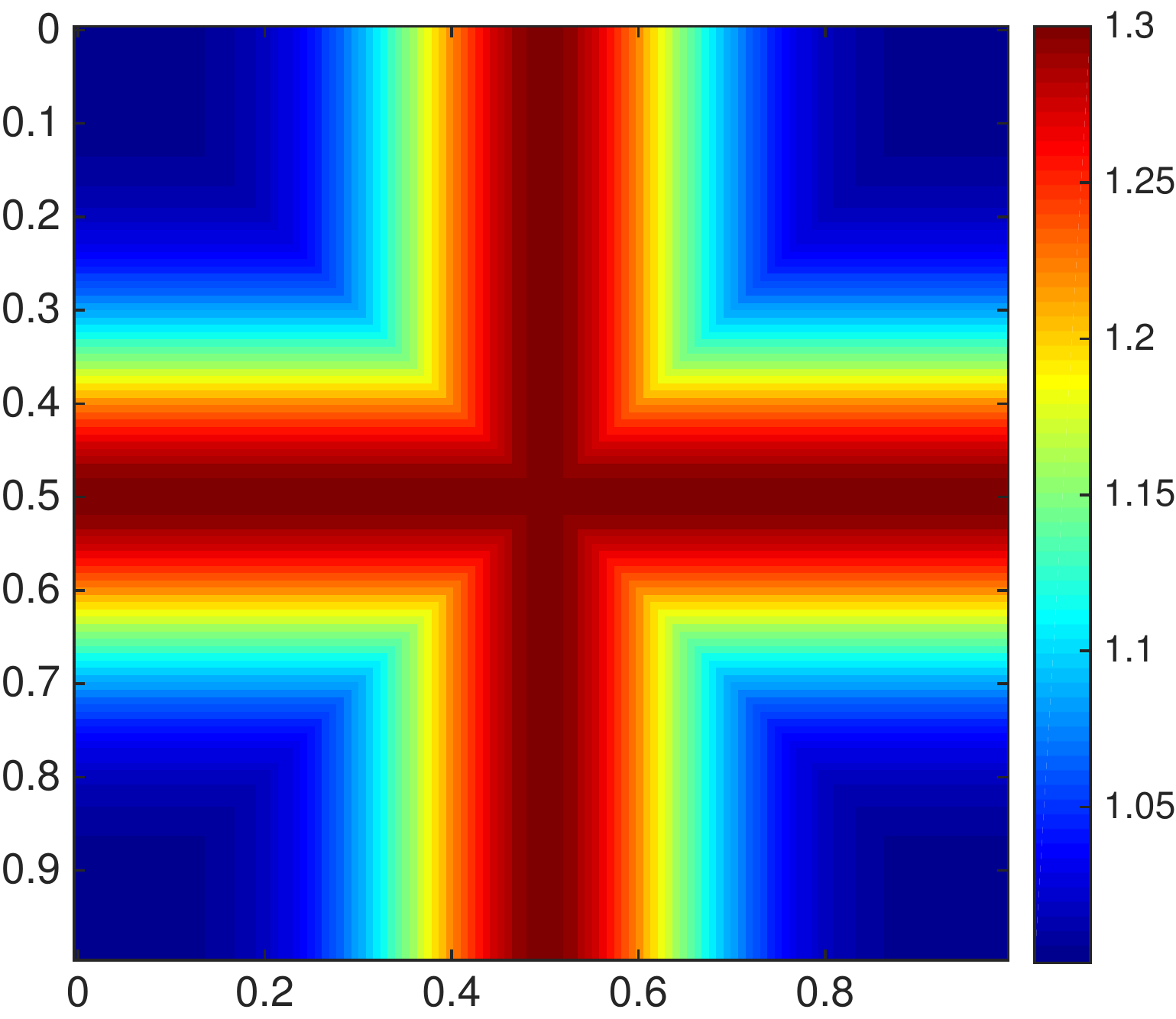}
  \includegraphics[width=0.425\textwidth]{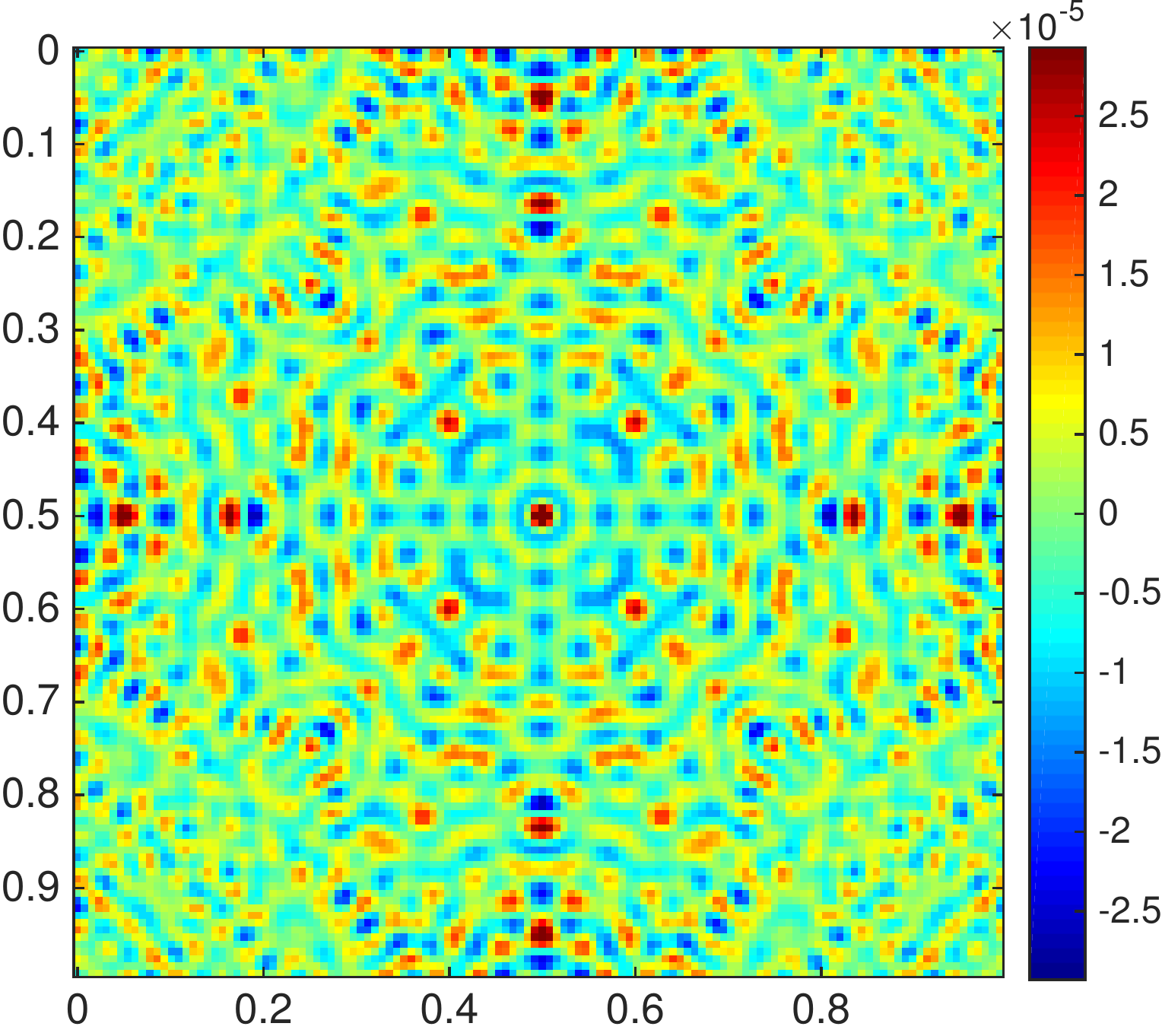}
  \caption{Results for velocity field (\ref{3DHii}) of the 3D Helmholtz equation. Top: numerical results. Bottom: $c(x)$ (left) and $u(x)$ (right) at $x_3 = 0.5$ for the largest problem size.}
  \label{tab:3DHii}
\end{table}

From the numerical tests for the Helmholtz equation one can make the
following two observations:
\begin{enumerate}
\item
  The iteration numbers are not sensitive to the growth of the problem
  size. In almost all cases, the preconditioned iterative solver
  converges in about 6-7 iterations. This clearly shows significant
  improvement over the results reported in \cite{spspd}.
\item
  The construction cost of the stencils is dominated by the setup cost
  of the nested dissection algorithm. This domination is more
  noteworthy in 3D due to a larger scaling difference between the
  stencil construction cost and the setup cost of the nested
  dissection factorization.
\end{enumerate}

\paragraph{Schr\"odinger Equation.}
For the Schr\"odinger equation, we set the system size to be
$n=1/h$. With the right rescaling, $v(x) = v_{\text{ext}}(x/h)/h^2 -
E/h^2$ where $v_{\text{ext}}(x)$ is the external potential field and
$E$ is the energy shift. We set $E=2.4$ so that there are at least
four points per oscillation. The potential fields tested for 2D are
\begin{enumerate}
[(i).]
\item
\label{2DSi}
An array of randomly put 2D Gaussians in the square.
\item
\label{2DSii}
An equal spaced array of 2D Gaussians with one missing at the center of the square.
\end{enumerate}
The results are given in Tables \ref{tab:2DSi} and \ref{tab:2DSii}.

\begin{table}[ht!]
  \centering
  \begin{tabular}{cc|cc|cc}
    \hline
    $N$ & $|S|$ & $T_{\text{stencil}}$ & $T_{\text{NDsetup}}$ & $N_{\text{iter}}$ & $T_{\text{NDsolve}}$ \\
    \hline
$64^2$ & $4$ & $2.11e-02$ & $5.89e-02$ & $7$ & $5.58e-02$ \\ 
$128^2$ & $8$ & $3.66e-02$ & $2.28e-01$ & $8$ & $2.28e-01$ \\ 
$256^2$ & $16$ & $9.70e-02$ & $9.76e-01$ & $8$ & $8.80e-01$ \\ 
$512^2$ & $32$ & $3.83e-01$ & $4.28e+00$ & $10$ & $4.18e+00$ \\ 
    \hline
  \end{tabular}
  \includegraphics[width=0.425\textwidth]{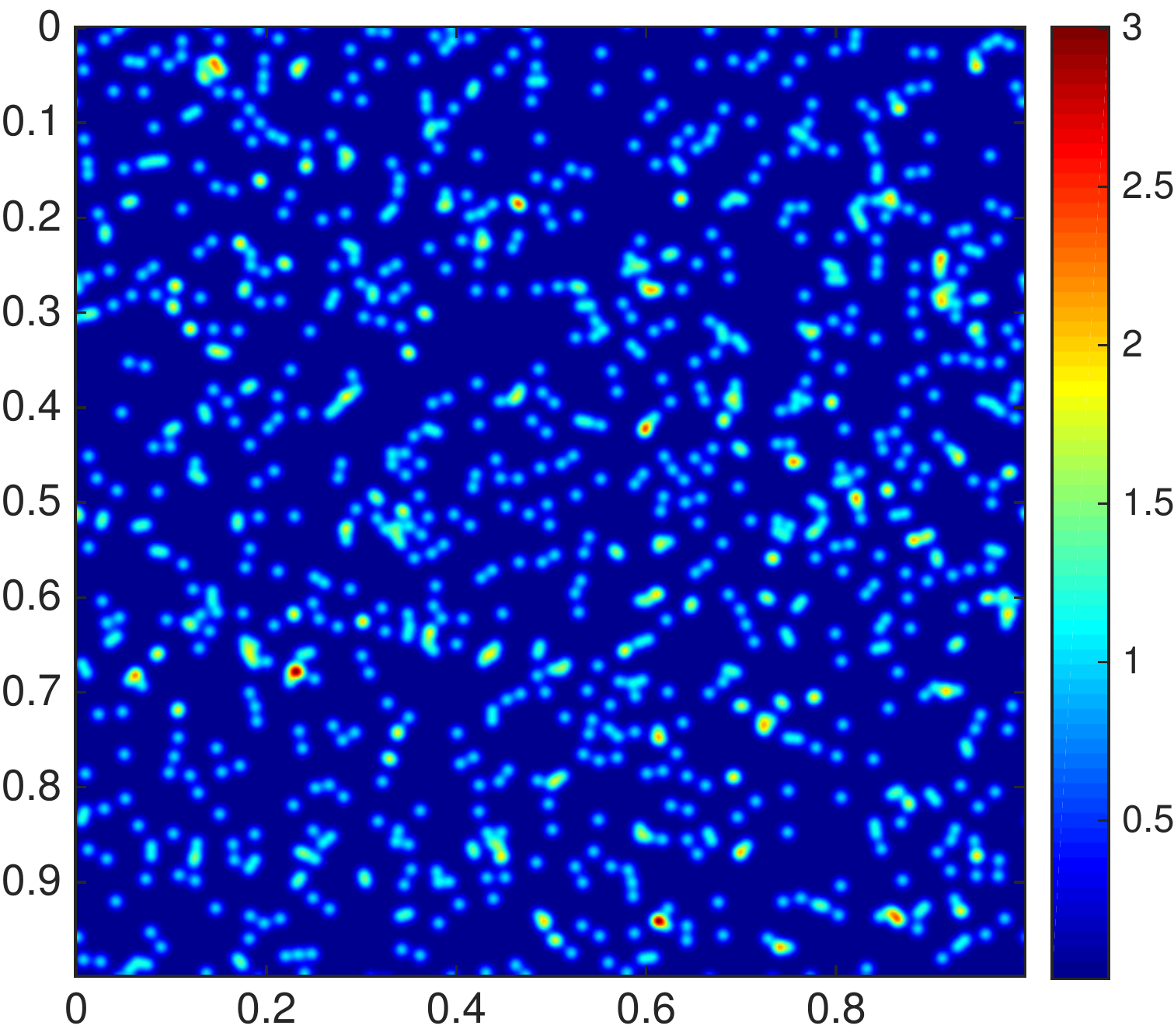}
  \includegraphics[width=0.425\textwidth]{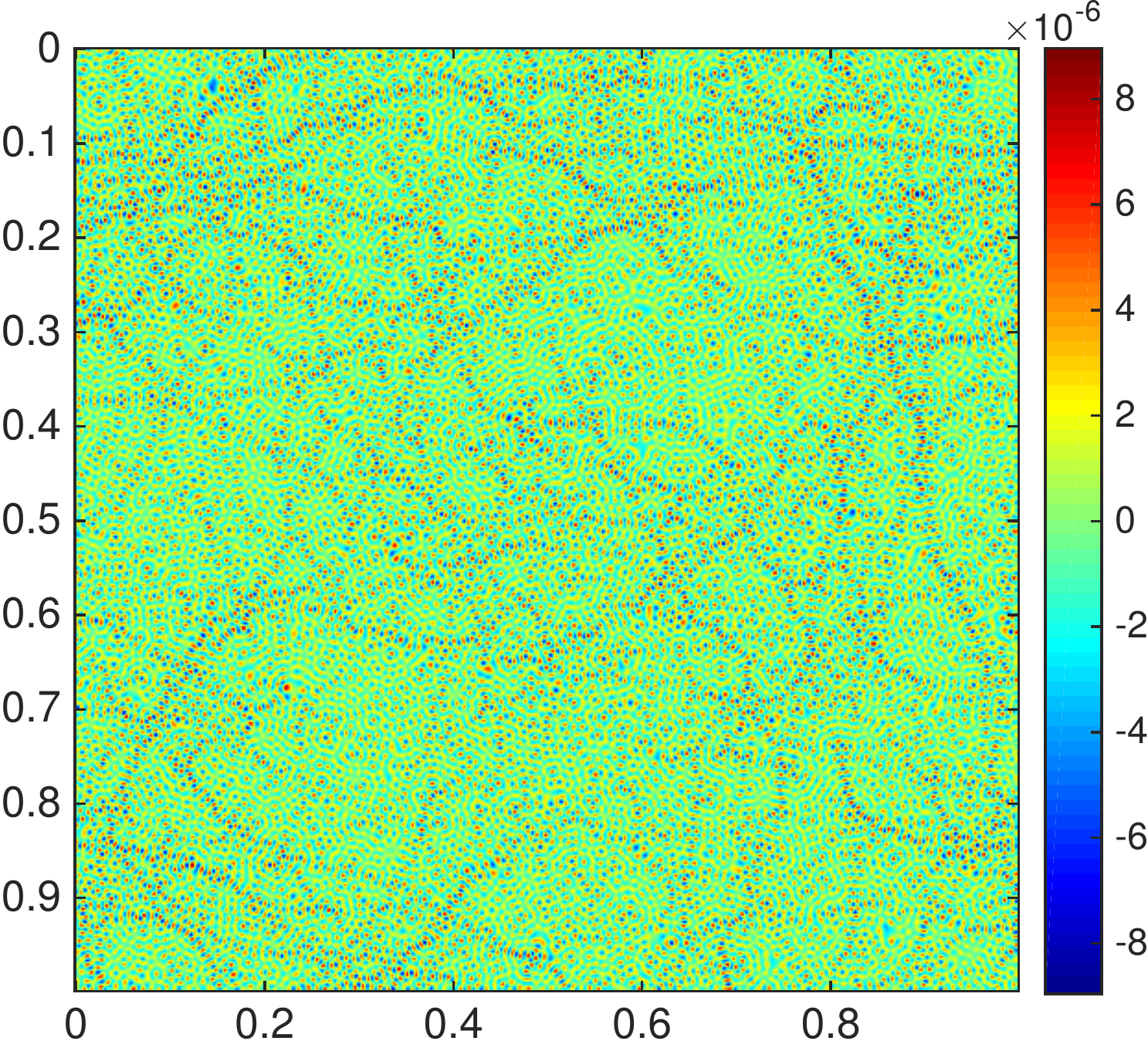}
  \caption{Results for potential field (\ref{2DSi}) of the 2D Schr\"odinger equation. Top: numerical results. Bottom: $v_{\text{ext}}(x/h)$ (left) and $u(x)$ (right) for the largest problem size.}
  \label{tab:2DSi}
\end{table}

\begin{table}[ht!]
  \centering
  \begin{tabular}{cc|cc|cc}
    \hline
    $N$ & $|S|$ & $T_{\text{stencil}}$ & $T_{\text{NDsetup}}$ & $N_{\text{iter}}$ & $T_{\text{NDsolve}}$ \\
    \hline
$64^2$ & $4$ & $1.89e-02$ & $5.35e-02$ & $6$ & $4.46e-02$ \\ 
$128^2$ & $8$ & $3.06e-02$ & $2.32e-01$ & $7$ & $2.12e-01$ \\ 
$256^2$ & $16$ & $9.47e-02$ & $8.59e-01$ & $7$ & $7.89e-01$ \\ 
$512^2$ & $32$ & $3.46e-01$ & $4.23e+00$ & $9$ & $3.69e+00$ \\ 
    \hline
  \end{tabular}
  \includegraphics[width=0.425\textwidth]{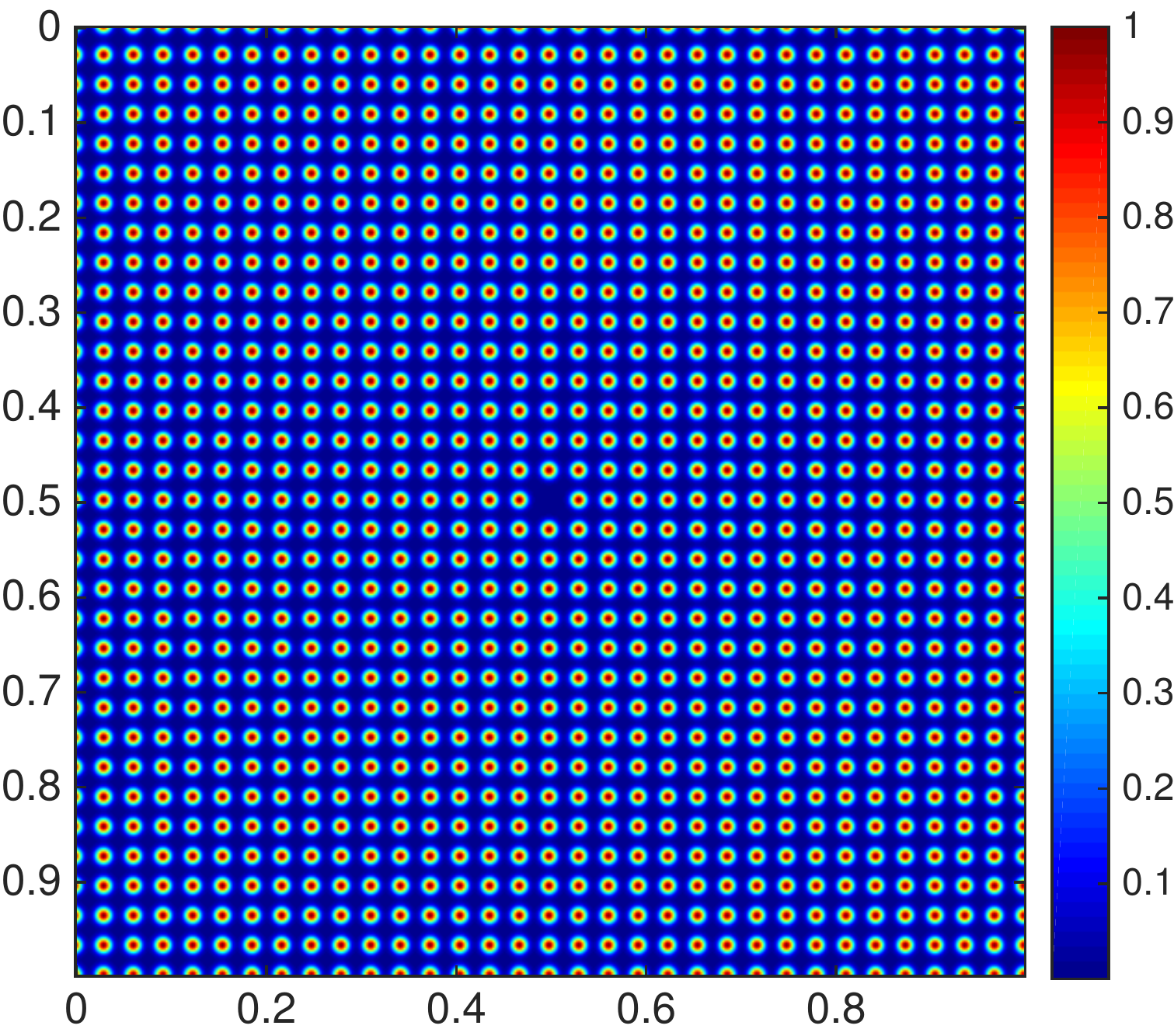}
  \includegraphics[width=0.425\textwidth]{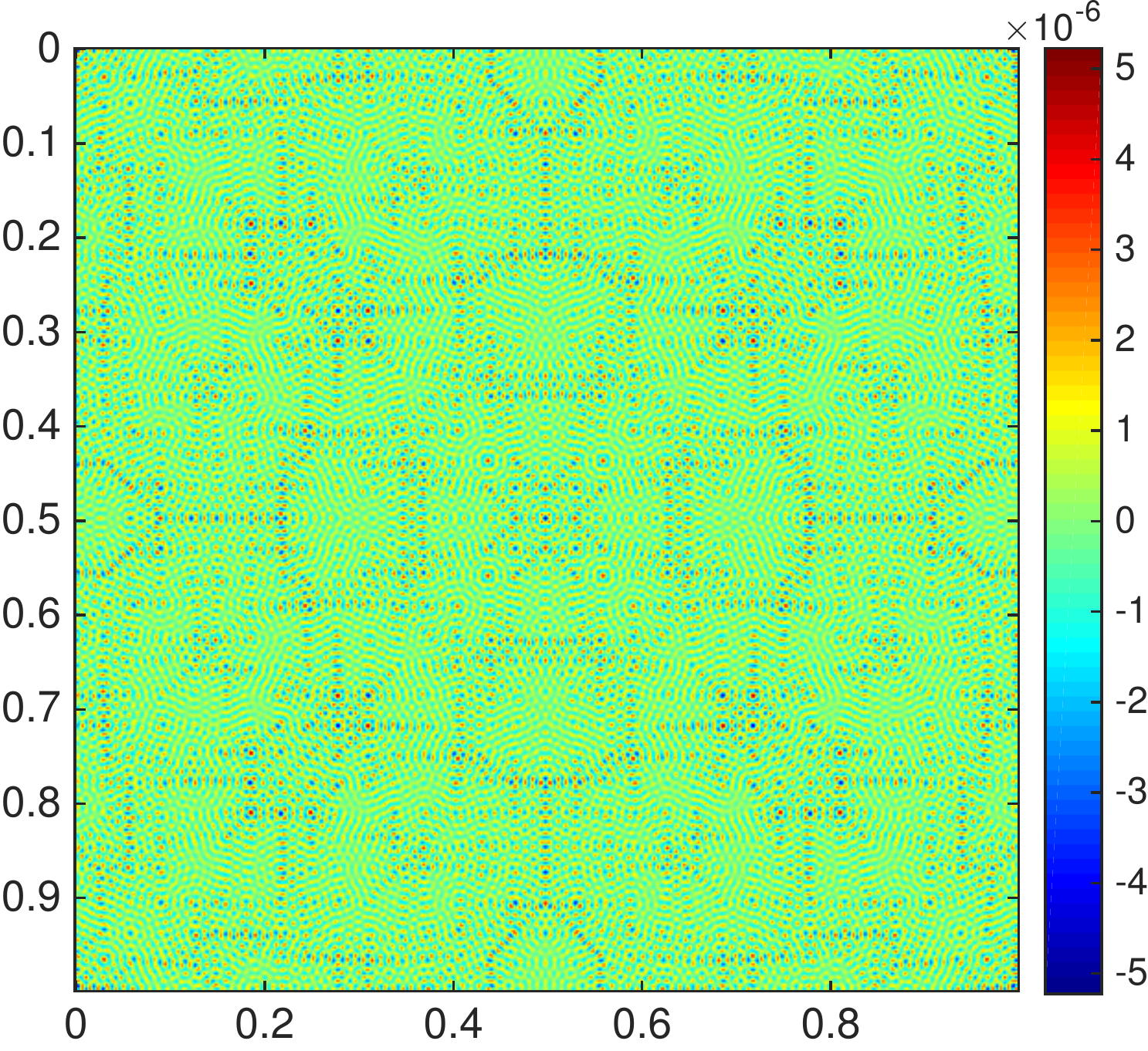}
  \caption{Results for potential field (\ref{2DSii}) of the 2D Schr\"odinger equation. Top: numerical results. Bottom: $v_{\text{ext}}(x/h)$ (left) and $u(x)$ (right) for the largest problem size.}
  \label{tab:2DSii}
\end{table}

For 3D tests, the potential fields are
\begin{enumerate}[(i).]
\item\label{3DSi} An array of randomly put 3D Gaussians in the cube.
\item\label{3DSii} An equal spaced array of 3D Gaussians with one
  missing at the center of the cube.
\end{enumerate}
The results are shown in Tables \ref{tab:3DSi} and \ref{tab:3DSii}.

\begin{table}[ht!]
  \centering
  \begin{tabular}{cc|cc|cc}
    \hline
    $N$ & $|S|$ & $T_{\text{stencil}}$ & $T_{\text{NDsetup}}$ & $N_{\text{iter}}$ & $T_{\text{NDsolve}}$ \\
    \hline
$16^3$ & $4$ & $3.53e-02$ & $2.58e-01$ & $7$ & $7.53e-02$ \\ 
$32^3$ & $8$ & $8.70e-02$ & $6.16e+00$ & $12$ & $1.86e+00$ \\ 
$64^3$ & $16$ & $3.99e-01$ & $1.68e+02$ & $9$ & $1.57e+01$ \\ 
$128^3$ & $32$ & $3.33e+00$ & $8.65e+03$ & $10$ & $2.19e+02$ \\ 
    \hline
  \end{tabular}
  \includegraphics[width=0.425\textwidth]{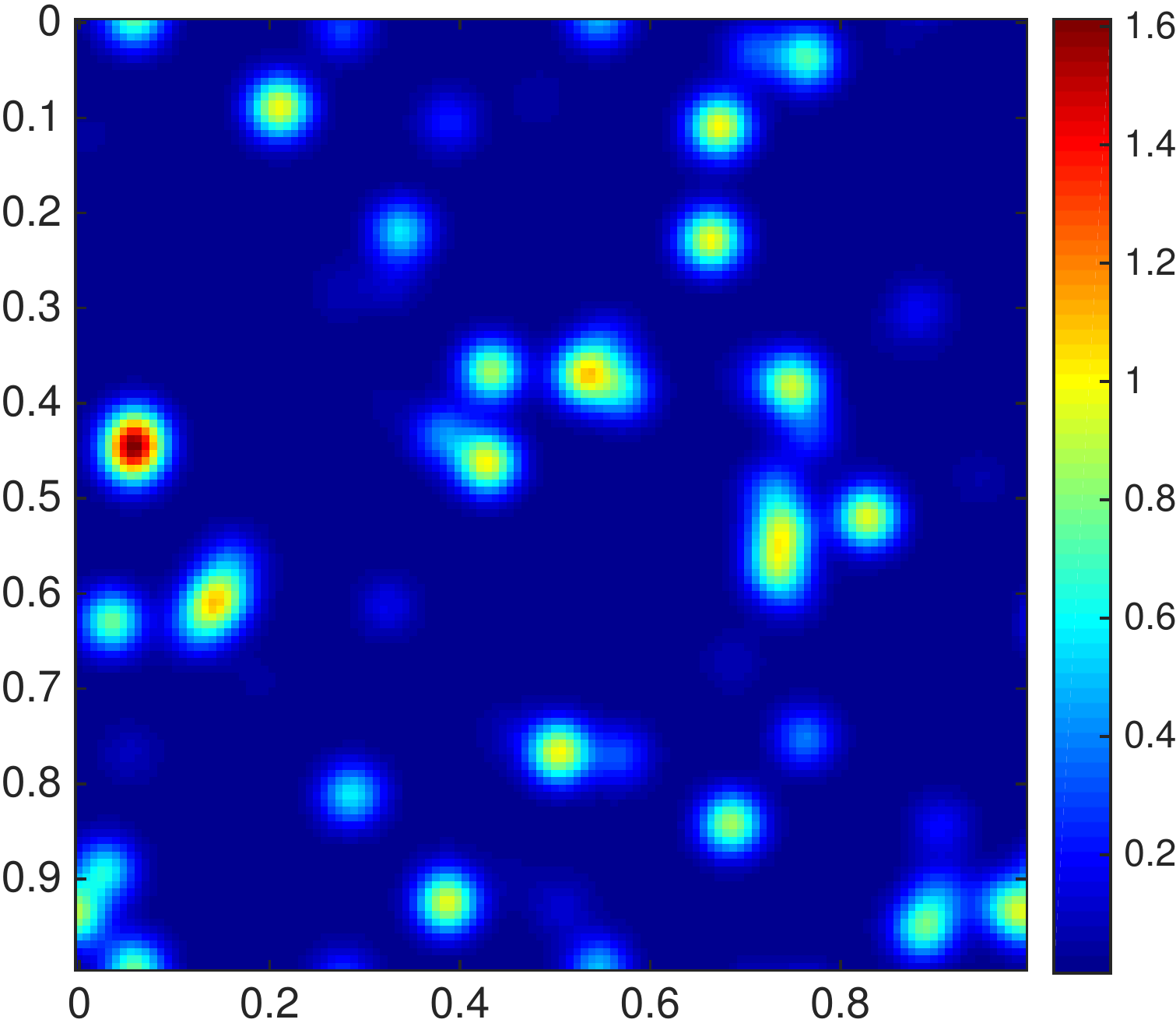}
  \includegraphics[width=0.425\textwidth]{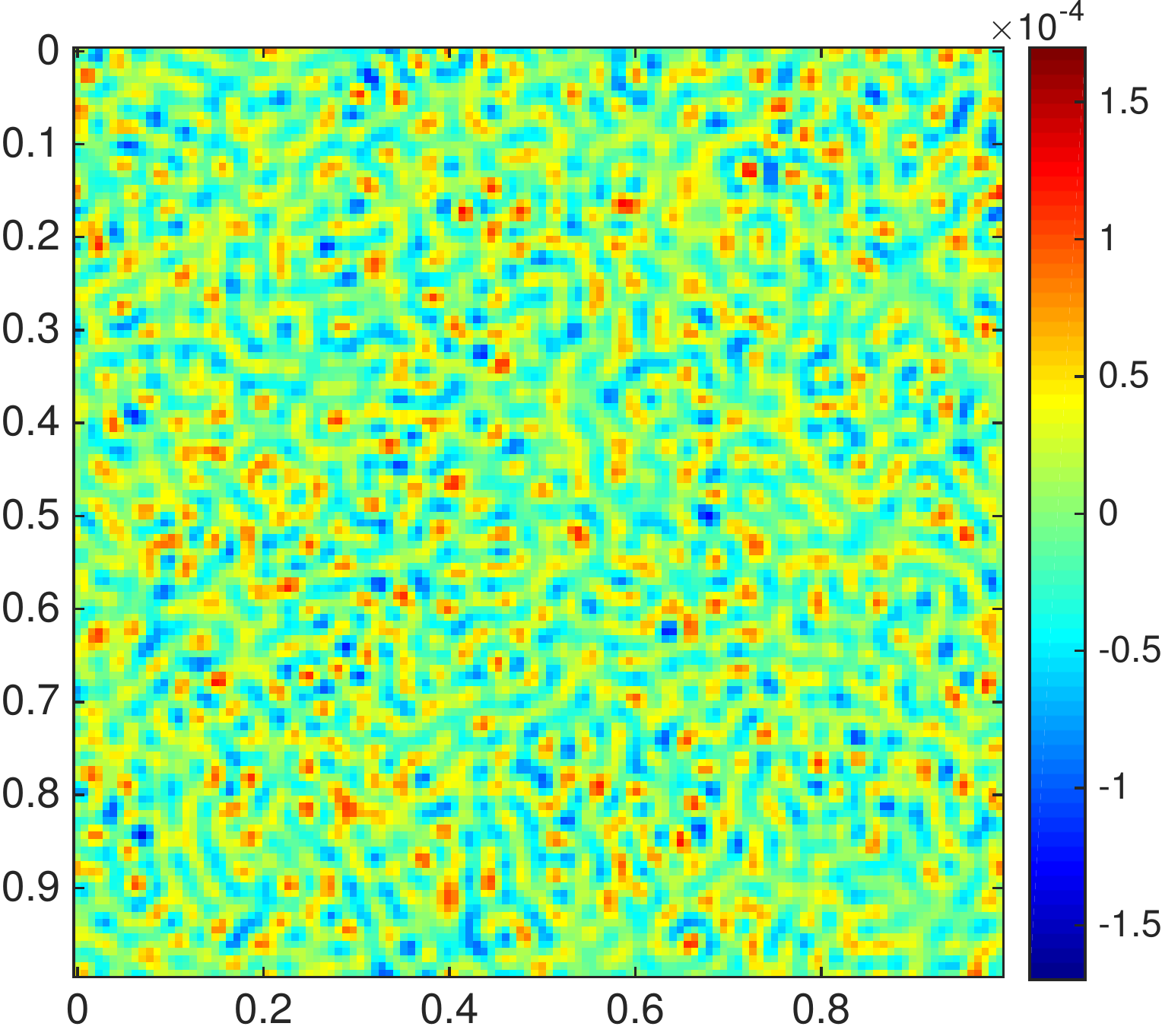}
  \caption{Results for potential field (\ref{3DSi}) of the 3D
    Schr\"odinger equation. Top: numerical results. Bottom:
    $v_{\text{ext}}(x/h)$ (left) and $u(x)$ (right) at $x_3 = 0.5$ for
    the largest problem size.}
  \label{tab:3DSi}
\end{table}

\begin{table}[ht!]
  \centering
  \begin{tabular}{cc|cc|cc}
    \hline
    $N$ & $|S|$ & $T_{\text{stencil}}$ & $T_{\text{NDsetup}}$ & $N_{\text{iter}}$ & $T_{\text{NDsolve}}$ \\
    \hline
$16^3$ & $4$ & $3.16e-02$ & $2.85e-01$ & $6$ & $7.16e-02$ \\ 
$32^3$ & $8$ & $1.08e-01$ & $6.14e+00$ & $7$ & $1.06e+00$ \\ 
$64^3$ & $16$ & $3.99e-01$ & $1.66e+02$ & $7$ & $1.32e+01$ \\ 
$128^3$ & $32$ & $3.61e+00$ & $8.54e+03$ & $7$ & $1.48e+02$ \\ 
    \hline
  \end{tabular}
  \includegraphics[width=0.425\textwidth]{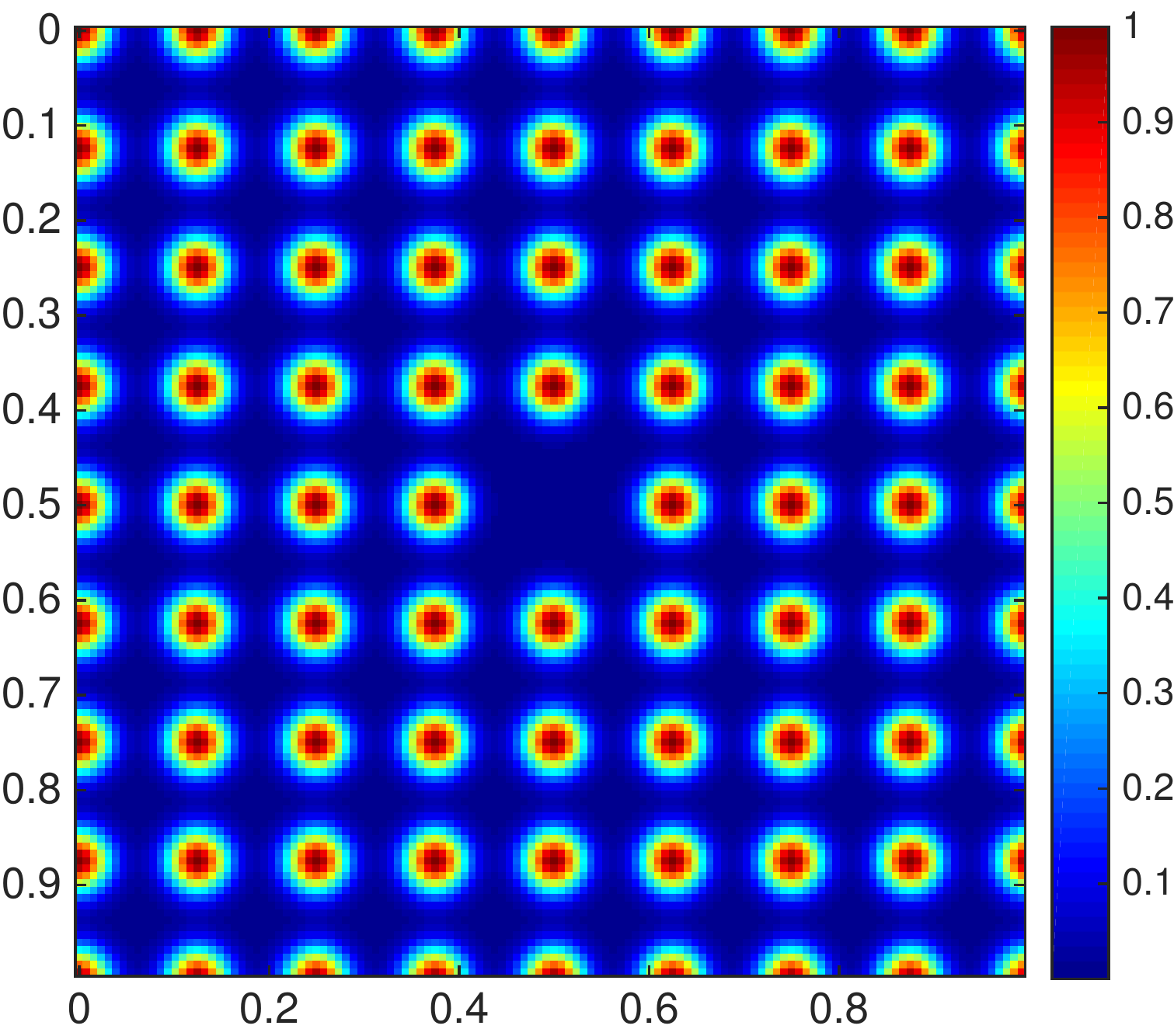}
  \includegraphics[width=0.425\textwidth]{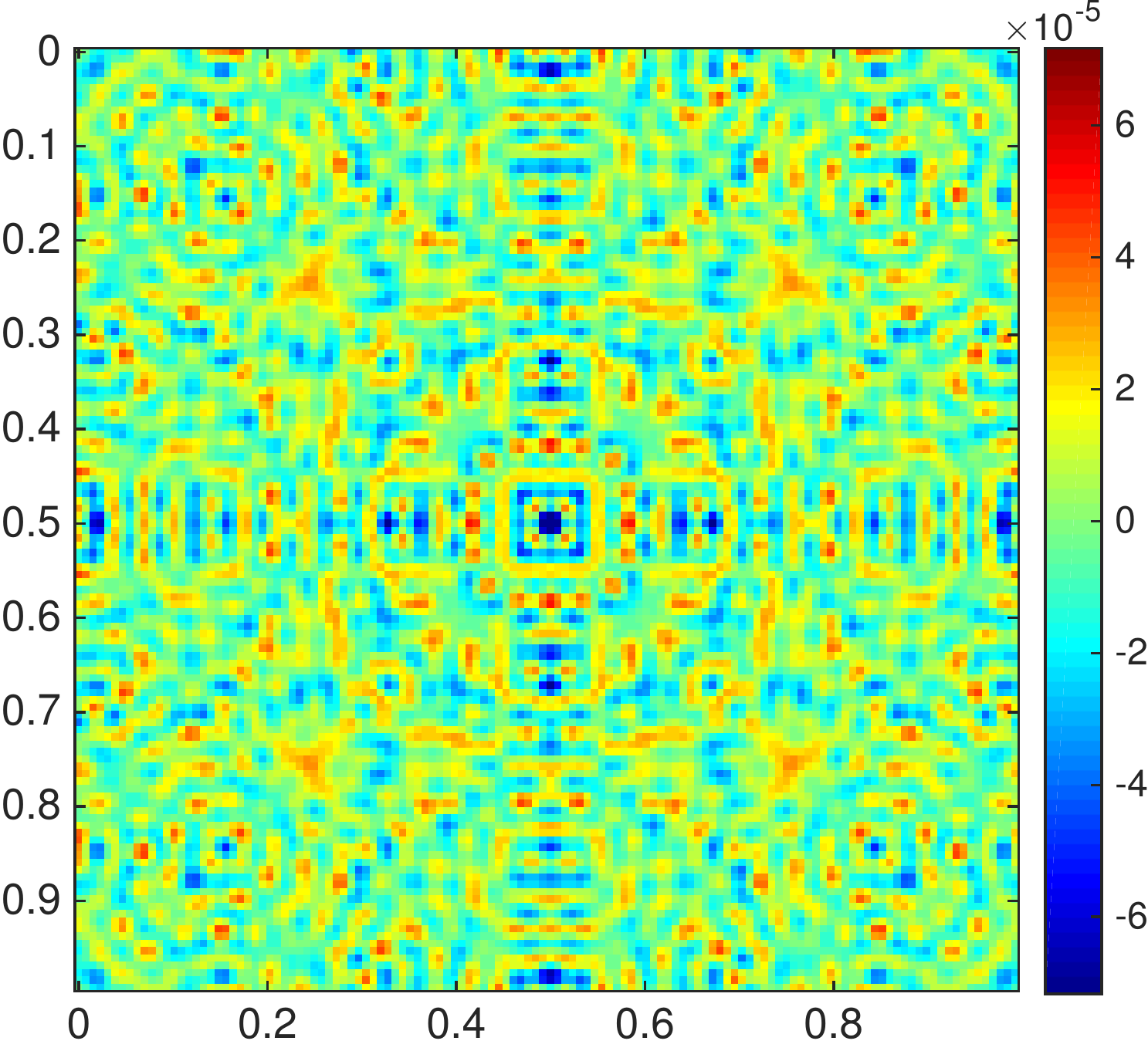}
  \caption{Results for potential field (\ref{3DSii}) of the 3D
    Schr\"odinger equation. Top: numerical results. Bottom:
    $v_{\text{ext}}(x/h)$ (left) and $u(x)$ (right) at $x_3 = 0.5$ for
    the largest problem size.}
  \label{tab:3DSii}
\end{table}

For the Schr\"odinger equation, the iteration number grows slightly
faster. The reason is that the variation of the field $v(x)$ increases
as the problem size increases in the Shcr\"odinger equation case,
while for the Helmholtz equation only $\omega$ grows with the problem
size and $c(x)$ remains the same. Thus the tests for the Schr\"odinger
equation are more challenging due to larger local variations of the potential field. Nevertheless, the growths of the
iteration numbers in the tests of the Schr\"odinger equation are still
mild compared to the growth of the problem size.

\section{Conclusion}
\label{sec:conclusion}
%
%
%

This paper introduces the localized sparsifying preconditioner for the
pseudospectral approximations of indefinite systems on periodic
structures based on the preconditioners in \cite{spspd,spspc}. The
novelty includes two parts. First, the local potential information is
taken into consideration during the construction of the sparse
matrices, which lowers the iteration number. Second, an FFT based
approach is introduced to compute the stencil which improves the
efficiency of the setup process of the algorithm.

Numerical tests show that the iteration number grows only mildly when
the problem size increases, which implies that solving pseudospectral
approximations of indefinite systems on periodic structures is not
inherently harder than solving sparse systems, up to a mildly growing
factor.

Another advantage of this new preconditioner is that, the construction
of the stencils of the algorithm is independent of the setup stage of
the nested dissection algorithm. The potential shift list $S$ needs
little information about the actual potential field $v(x)$ except for
the minimum and the maximum value, which means that the stencils can
be built in advance, and as long as the value of $v(x)$ is in a
certain range, there is no need to reconstruct the stencil no matter
how $v(x)$ varies. This can be helpful when an iterative process is
involved or $v(x)$ is constantly changing, such as in
\cite{lu2015sparsifying}.

The choice of $t$ for each location in the current setting is rather
crude. There are several ways to make potential improvements. For
example, one can adopt the stencil where $t=1$ only for the largest
skeleton in the nested dissection algorithm, while for the rest of the
points, stencils with higher $t$ values can be used. In this way, the
setup cost of the nested dissection algorithm will not increase too
much, while the iteration number may be further reduced. One can also
use the stencils with lower $t$ values for the locations where the
local potential field variation is milder and with higher $t$ values
where the variation is more drastic. These techniques can be helpful
for practical application of this algorithm.

\section*{Acknowledgments}
The authors are partially supported by the National Science
Foundation under award DMS-1521830 and the U.S. Department of Energy's
Advanced Scientific Computing Research program under award
DE-FC02-13ER26134/DE-SC0009409.

\bibliographystyle{abbrv}
\bibliography{references}
\end{document}